\documentclass{amsart} 

\usepackage{amsmath,amsthm,amssymb,amscd,enumitem,mathtools}
\usepackage{amsfonts}
\usepackage{rotating}
\usepackage{euscript}
\usepackage{pst-node}
\usepackage{epsfig,verbatim}
\mathtoolsset{showonlyrefs}

\def\be{\begin{equation}}
\def\ee{\end{equation}}

\def\C{{\mathbb C}} 
\def\f{\EuScript}
 
\def\P{{\mathbb P}}
\def\Z{{\mathbb Z}}

\def\phi{{\varphi}}
 
\def\tt{\widetilde}
\def\deg{{\rm deg\,}}

\def\GCD{{\rm GCD }}
\def\LCM{{\rm LCM }}

\def\bp{\begin{proposition}}
\def\ep{\end{proposition}}

\def\bt{\begin{theorem}}
\def\et{\end{theorem}}
\def\br{\begin{remark}}
\def\er{\end{remark}}
\def\be{\begin{equation}}

\def\bee{\begin{equation*}}
\def\l{\label}

\def\m{\mu}
\def\ee{\end{equation}}
\def\eee{\end{equation*}}
\def\bl{\begin{lemma}}
\def\el{\end{lemma}}
\def\bc{\begin{corollary}}
\def\ec{\end{corollary}}
\def\pr{\noindent{\it Proof. }}

\def\bd{\begin{definition}}
\def\ed{\end{definition}}
\def\t{\widetilde}

\def\Aut{{\rm Aut}}
\def\Mon{{\rm Mon}}

\newtheorem{theorem}{Theorem}[section]
\newtheorem{lemma}[theorem]{Lemma}
\newtheorem{definition}[theorem]{Definition}
\newtheorem{corollary}[theorem]{Corollary}
\newtheorem{proposition}[theorem]{Proposition}

\theoremstyle{definition}

\theoremstyle{definition}
\newtheorem{remark}[theorem]{Remark}

\begin{document}

\title[On iterates of rational functions]{On iterates of rational functions with maximal number of critical values}
\author{Fedor Pakovich}
\thanks{
This research was supported by ISF Grant  No. 1092/22}
\address{Department of Mathematics, Ben Gurion University of the Negev, Israel}
\email{
pakovich@math.bgu.ac.il}

\begin{abstract} Let $F$ be a rational function of one complex variable  of degree \linebreak $m\geq 2$.
The function $F$ is called simple if for every $z\in \C\P^1$ the preimage $F^{-1}\{z\}$ contains at least $m-1$ points.  We show that if $F$ is a simple rational function of degree $m\geq 4$ and 
$F^{\circ l} =G_r\circ G_{r-1}\circ \dots \circ G_1$, $l\geq 1$,
is a decomposition of an iterate  of $F$ into a composition of indecomposable  rational functions, then $r=l$ and there exist M\"obius transformations  $\mu_i,$ $1\leq i \leq r-1,$ such that 
$G_r=F\circ \mu_{r-1},$ $G_i=\mu_{i}^{-1}\circ F \circ \mu_{i-1},$ $1<i< r,$ and $G_1=\mu_{1}^{-1}\circ F$.  As applications, we solve a number of    
 problems in complex and arithmetic 
dynamics for ``general'' rational functions. 

\end{abstract}

\maketitle

\section{Introduction}

Let $F$ be a rational function of one complex variable of degree $m\geq 2.$ The function $F$ is called {\it indecomposable} if the equality $F=F_2\circ F_1$, where $F_1,$ $F_2,$ are rational functions,  implies that at least one of the functions $F_1,F_2$ is of degree one. 
Any representation of  $F$  in the form $F=F_r\circ F_{r-1}\circ \dots \circ F_1,$
where $F_1,F_2,\dots, F_r$ are rational functions of degree at least two,
is called {\it a decomposition} of $F.$ Two decompositions 
\be \l{lll} F=F_r\circ F_{r-1}\circ \dots \circ F_1 \ \ \ \ {\rm and} \ \ \ \  
F=G_{l}\circ G_{l-1}\circ \dots \circ G_1\ee 
are called {\it equivalent} if $l=r$ and either $r=1$ and $F_1=G_1$, or $r\geq 2$ and there exist M\"obius transformations  $\mu_i,$ $1\leq i \leq r-1,$ such that 
$$F_r=G_r\circ \mu_{r-1}, \ \ \ 
F_i=\mu_{i}^{-1}\circ G_i \circ \mu_{i-1}, \ \ \ 1<i< r, \ \ \ {\rm and} \ \ \ F_1=\mu_{1}^{-1}\circ G_1.
$$ 
It is obvious that any 
rational function $F$ of degree $m\geq 2$ can be decomposed into a composition of indecomposable rational functions, although in general not in a unique way. 
The problem of describing all such decompositions is quite delicate, and the general theory exists only if $F$ is a polynomial or a Laurent polynomial (see \cite{r1}, \cite{pak}). 

In dynamical applications, one needs to have a description 
of decompositions of the whole totality of iterates of a given rational function $F$ (see e.g. \cite{cm}, \cite{gtz}, \cite{gtz2}, \cite{ms}, \cite{pj}, \cite{ic}, \cite{tame}), and  the main result of this paper states roughly speaking that for a general rational function of degree $m\geq 4$ all such decompositions are trivial. As applications, we solve a number of    
 problems in complex and arithmetic dynamics for general rational functions. 
Here and below, saying that some statement holds for {\it general} rational functions of degree $m$, we mean the following: if we  
identify the set of rational functions  of degree $m$ with an algebraic variety ${\rm Rat}_m$ obtained from $\C\P^{2m+1}$ by removing the resultant hypersurface,   then this statement holds for all $F\in {\rm Rat}_m$ with exception of some proper Zariski closed subset.

In more detail, we prove a number of results, which hold for {\it simple} rational functions, that is, for rational functions $F$ of degree $m\geq 2$ such that for every $z\in \C\P^1$ the preimage $F^{-1}\{z\}$ contains at least $m-1$ points. 

Our main result is the following statement.

\bt \l{t1} Let $F$ be a simple rational function of degree $m\geq  4$.
Then any decomposition of $F^{\circ l}$, $l\geq 1$, into a composition of indecomposable  rational functions is equivalent to $F^{\circ l}$.
\et

We apply Theorem \ref{t1} to describing a variety of objects associated with a rational function $F$  of degree at least two, using the following notation.
\vskip 0.2cm 
\begin{tabular}{@{ }ll}

$\langle F\rangle $ & is the semigroup of  rational functions generated by $F$. \\ 
    $C(F)$     & is the semigroup of all rational functions  commuting with $F$.\\
$\Aut(F)$ & is the group of all M\"obius transformations belonging to $C(F)$. \\
   $C_{\infty}(F)$ & is the semigroup of all rational functions commuting with \\ &  some iterate  
 of $F$.\\
$\Aut_{\infty}(F)$ & is the group of all M\"obius transformations belonging \\ & to $C_{\infty}(F)$. \\
$\langle \Aut_{\infty}(F), F \rangle$ & is the semigroup of  rational functions generated by $F$ \ \\ & and $\Aut_{\infty}(F).$ \\ 
$\mu_F$ & is the measure of maximal entropy of $F$. \\
$E_0(F)$ & is the group of all 
M\"obius transformations preserving $\mu_F$.\\ 
$E(F)$ & is the semigroup consisting of all rational functions $G$ of degree \\ & at   least two  with 
 $\mu_G=\mu_F$, completed by the group $E_0(F)$.\\ 
 $G_0(F)$ & is the maximal subgroup of $\Aut(\C\P^1)$ such that 
for every \\ & $\sigma\in G_0(F)$ 
there exists $\nu\in G_0(F)$ satisfying  
$ F\circ \sigma=\nu \circ F$.

\end{tabular}

\vskip 0.2cm

Using Theorem \ref{t1}, we show that for simple rational functions  the above objects 
 are related in a very simple way. 

\bt \l{t2} Let $F$ be a simple rational function of degree $m\geq 4 $. Then 
$$E_0(F)=\Aut_{\infty}(F)=G_0(F)\ \ \ \ {\it and} \ \ \ \ E(F)=C_{\infty}(F)=\langle \Aut_{\infty}(F), F \rangle. $$
\et

The link between Theorem \ref{t1} and Theorem \ref{t2} is based on the  results of Ritt (\cite{r}) and   Levin and  Przytycki (\cite{lev}, \cite{lp}). Namely, the theorem of Ritt about commuting rational functions implies that for a fixed non-special rational function $F$ of degree at least two,  a rational function $G$ of degree at least two belongs to  
$C_{\infty}(F)$ if and only if the equality 
\be \l{u} F^{\circ k}=G^{\circ l}\ee 
holds for some $k,l\geq 1$. 
On the other hand, the  results of   Levin and  Przytycki yield that   
$G$ belongs to  $E(F)$ if and only if the equality 
\be \l{uu} F^{\circ k_1}=F^{\circ k_2}\circ G^{\circ l} \ee 
holds for some $k_1,l\geq 1$, $k_2\geq 0$ (see Section \ref{rpr} for more detail).

As an application of  Theorem \ref{t2}, we prove the following result.

\bt \l{t3} For a general rational function $F$  of degree $m\geq  4$, the equalities  
$$E_0(F)=\Aut_{\infty}(F)=G_0(F)=id\ \ \ \ {\it and} \ \ \ \ E(F)=C_{\infty}(F)=\langle  F \rangle$$ hold. 
\et

Notice that Theorem \ref{t3} 
provides an affirmative answer to the question of Ye, who proved that the equality $E(F)=\langle F\rangle$ holds after removing  from ${\rm Rat}_m$  countably many algebraic sets, and asked whether it remains true if to remove from ${\rm Rat}_m$ only finitely many such sets (\cite{ye}).

Further applications of Theorem \ref{t1} concern problems that can be reformulated in terms of semiconjugacies between rational functions (see the papers 
\cite{bu}, \cite{e3}, \cite{fg}, \cite{i},  \cite{ms}, \cite{ac}, \cite{ic} for examples of such problems).  We recall that 
a rational function $B$ of degree at least two is called {\it semiconjugate} to a rational  function $A$
if there exists a non-constant rational function $X$
such that the diagram 
\be \l{ii1} 
\begin{CD}
\C\P^1 @>B>> \C\P^1 \\
@V X VV @VV  X V\\ 
\f \C\P^1 @>A>> \f\C\P^1\
\end{CD}
\ee
commutes. 
 A comprehensive description of triples $A,B,X$ such that \eqref{ii1} commutes was obtained in the series 
of papers \cite{semi}, \cite{rec}, \cite{fin}, \cite{lattes}. For simple $A$, Theorem \ref{t1} permits to reduce this description  
to the following uncomplicated form suitable for applications. 


\bt \l{t4} Let $F$ be a simple rational function of degree $m\geq 4$, and 
$G,X$ non-constant rational functions  such that  the diagram 
\be \l{ii2}
\begin{CD}
\C\P^1 @>G>> \C\P^1 \\
@V X VV @VV  X V\\ 
\f \C\P^1 @>F^{\circ r}>> \f\C\P^1\
\end{CD}
\ee
commutes for an integer $r\geq 1$. Then there exist a M\"obius transformation $\nu$ and an integer $l\geq 0$ such that the equalities 
$$X=F^{\circ l}\circ \mu, \ \ \ \ \ \ \ G=\mu^{-1} \circ F^{\circ r} \circ \mu$$ 
hold.
\et

As an example of an application of Theorem \ref{t4}, we consider the problem of describing  
periodic algebraic curves for endomorphisms  of $(\C\P^1)^2$ of the form
\be \l{ss} (F_1,F_2):\, (z_1,z_2)\rightarrow (F_1(z_1),F_2(z_2)),\ee
where  $F_1,F_2$ are rational functions, which   reduces to describing solutions of a system of semiconjugacies.  A description of periodic curves in the case  where  $F_1,F_2$ are polynomials was obtained by Medvedev and Scanlon (\cite{ms}) and has numerous applications in complex and arithmetic dynamics 
(see e. g. \cite{bm}, \cite{fg},  \cite{gn}, \cite{gn1},  \cite{gn2}, \cite{gx}, \cite{n}).
A description of periodic curves in the general case     
was obtained in the recent paper \cite{ic}. Notice that the problem of describing periodic curves is closely related to 
 a variant of a conjecture of Zhang (\cite{z}) 
on the existence of Zariski dense orbits for endomorphisms $(F_1,F_2)$ defined over a field $K$ of characteristic zero (see \cite{ms}, \cite{ic}, \cite{k}).

Theorem \ref{t4} permits to shorten considerably the results of \cite{ic}  in case $F_1$ and $F_2$ are simple, leading to the following result, which can be easily used  for applications.

\bt \l{t5}  Let $F_1$ and $F_2$ be simple rational functions of degree $m\geq 4$, and 
$C$ an irreducible algebraic  curve in $(\P^1(\C))^2$  that is not a vertical or horizontal line. Then $(F_1,F_2)^{\circ d}(C)=C$ for an integer $d\geq 1$  if and only if 
\be \l{moca} F_2^{\circ d}=\alpha \circ F_1^{\circ d}\circ\alpha^{-1}\ee for some  M\"obius transformation 
$\alpha$, and $C$ is one of the graphs $$y=(\alpha\circ \mu\circ F_1^{\circ s})(x), \ \ \ \ \ x=(\mu\circ F_1^{\circ s}\circ\alpha^{-1})(y),$$  where $\mu\in \Aut(F_1^{\circ d})$ and $s\geq 0.$ 
\et

Using Theorem \ref{t5}, we prove the following result about invariant and periodic curves for 
general rational functions. 

\bt \l{t6} For every $m\geq 4$ there exists a Zariski open set $U$ in  ${\rm Rat}_m$ such that the following holds. For any $F_1,F_2\in U$, an irreducible algebraic  curve $C$ in $(\P^1(\C))^2$  that is not a vertical or horizontal line   is $(F_1,F_2)$-periodic if and only if 
$$F_2=\alpha \circ F_1\circ\alpha^{-1}$$ for some  M\"obius transformation 
$\alpha$, and $C$ is one of the graphs $$y=(\alpha\circ  F_1^{\circ s})(x), \ \ \ \ \ x=(F_1^{\circ s}\circ\alpha^{-1})(y),$$  where  $s\geq 0.$ In particular, any $(F_1,F_2)$-periodic curve is $(F_1,F_2)$-invariant. 
\et

For proving Theorem \ref{t1}, we use the following strategy.  
First, we show that if $F$ is a simple rational function of degree $m\geq 4$ and 
$H$ is an indecomposable rational function of degree at least two such that the algebraic curve 
\be \l{cura} H(y)-F(x)=0\ee 
is irreducible, then the genus of this curve is greater than zero. Second, we show that if  \eqref{cura} is reducible, then either $H=F\circ \mu,$ where $\mu$ is a M\"obius transformation, or  $\deg H$ is equal to the binomial coefficient $\binom{m}{k}$ for some $k,$ $1<k <m-1.$ Third, using the theorem of 
Sylvester \cite{syl} and Schur \cite{sch} about prime divisors of binomial coefficients, we show that  there exists a prime number $p$ such that $p\mid \binom{m}{k}$ but $p\nmid m.$ The above statements yield 
that if $ F^{\circ l}=H\circ R$ for some rational function $R$, then 
$H$ necessarily has the form $H=F\circ \mu$ for some M\"obius transformation $\mu,$  
and this fact allows us to prove the theorem.

The paper is organized as follows. In the second section, using the above approach  we prove Theorem \ref{t1}. 
In the third section, we deduce from Theorem \ref{t1}  Theorem \ref{t2} and Theorem \ref{t3}. In the fourth section, basing on results about semiconjugate rational functions and invariant curves from \cite{semi}, \cite{lattes}, \cite{ic}, we  prove Theorem \ref{t4},  Theorem \ref{t5},  and Theorem \ref{t6}.

Finally, in the fifth section, we 
give a number of conditions implying that some iterate $F^{\circ k}$, $k>1,$ of an indecomposable rational function $F$ has  a decomposition not equivalent to $F^{\circ k}$ itself.  
We also construct explicit examples of simple rational functions of degree 2 and 3 for 
which Theorems \ref{t1} - \ref{t2} and Theorems \ref{t4} - \ref{t5} are not true. As for Theorem \ref{t3} and Theorem \ref{t6}, we believe that they have some analogues for $m=2$ and $m=3$. However, the methods of this paper do not apply to this situation. 

\section{Decompositions of iterates of rational functions} 

\subsection{The monodromy group and decompositions} 
Let $G$ be a group which acts transitively on a finite set $S$. 
We recall that a subset $T$ of $S$ is called a block of  $G$, if for each $g\in G$ either $g(T)=T$ or $g(T)\cap T=\emptyset.$  
Clearly, if $T$ is a block, then $\f T=\{\sigma(T), \sigma\in G\}$ is a partition of $S$, which is called an imprimitvity system of $G$.  The group $G$ is called primitive if its blocks are only singletons and the whole $S$. Otherwise, $G$ is called imprimitive.

Let $F$ be a rational function, and 
$c(F)=\{z_1, z_2, \dots , z_r\}$ the set of all critical values of $F$. Let us fix a point 
$z_0\in\C\P^1\setminus c(F)$ and some loops $\gamma_i$ around $z_i,$ $1\leq i \leq r,$ such that 
$\gamma_1\gamma_2...\gamma_r=1$ in $\pi_1(\C\P^1\setminus c(F),z_0)$. Further, let us   
 denote by 
$\delta_i,$ $1\leq i \leq r,$ a permutation of points of 
$F^{-1}\{z_0\}$ induced by the lifting of $\gamma_i,$ $1\leq i \leq r$. In this notation, 
the monodromy group of $F$ is  defined as 
the permutation group  generated by $\delta_i,$ $1\leq i \leq r.$ We will  denote this group by $\Mon(F)$.

The imprimitivity systems of the group $\Mon(F)$ correspond to decompositions of $F$.  
Namely, if $F=A\circ B$ is a decomposition of $F$ into a composition of rational functions $A$ and $B$, where $\deg A=d,$ 
then $\Mon(F)$ has an imprimitivity system consisting of 
$d$ blocks $B^{-1}\{t_i\},$ $1\leq i \leq d,$ where $\{t_1,t_2,\dots, t_{d}\} =A^{-1}\{z_0\}$. 
Furthermore, any imprimitivity system of $\Mon(F)$ arises from a decomposition of $F$,
and to decompositions $F=A\circ B$ and $F=C\circ D$ corresponds the same
imprimitivity system if and only there exists a M\"obius transformation $\mu$ 
such that $$ A=C \circ \mu^{-1}, \ \ \ B=\mu \circ D.$$
In particular, $F$ is indecomposable if and only if $\Mon(F)$ is primitive.

\subsection{A calculation of the genus of $H(x)-F(y)=0$}  
Let $F$ be a rational function of degree $m\geq 2$. We denote 
by $\deg_zF$ the multiplicity of $F$ at a point $z\in \C\P^1.$   The following two results are known. We include the proofs for the reader convenience. 

\bl \l{bl} Let $F$ be a rational function of degree $m\geq 2$. Then the following conditions are equivalent. 
\begin{enumerate}[label=\roman*)]
\item The function is simple. 
\item The number of critical points of $F$ is equal to the number of critical values, and the multiplicity of $F$ at every critical point is equal to two. 
\item The number of critical values of $F$ is equal to $2m-2.$
\end{enumerate} 
\el
\pr 
The equivalence $i)\Leftrightarrow ii)$ follows from the definition. 
Furthermore, it follows from the Riemann-Hurwitz formula 
\be \l{rh} 2m-2=\sum_{z\in \C\P^1}(\deg_zF-1)\ee
 that the number of critical points of $F$  does not exceed $2m-2$, and the equality is attained if and only if the multiplicity of $F$ at every critical point is equal to two. Since 
 the number of critical values of $F$ does not exceed the number of critical points, this 
 implies easily the equivalence $ii)\Leftrightarrow iii).$
 \qed

\bt \l{0} Let $F$ be a simple  rational function of degree $m\geq 2$. Then $F$ is indecomposable, and  $\Mon(F)\cong S_m.$ 
\et 
\pr Assume that \be \l{bcr} F=F_1\circ F_2,\ee where $F_1$ and $F_2$ are rational functions of degrees $m_1$ and $m_2.$ Since  $F$ is simple, the number of critical values of $F$ is $2m_1m_2-2$ by Lemma \ref{bl}.  On the other hand, it follows from \eqref{bcr} by the chain rule that the  number of critical values of $F$ does not exceed $(2m_1-2)+(2m_2-2)$. Thus, 
$$2m_1m_2-2\leq (2m_1-2)+(2m_2-2),$$ implying that 
$$2m_1m_2-2-(2m_1-2)-(2m_2-2)=2(m_1-1)(m_2-1)\leq 0.$$ Therefore, at least one of the functions $F_1$ and $F_2$ has degree one. 

Since $F$ is indecomposable, the monodromy group $\Mon(F)$ of $F$ is primitive. Furthermore,  for any critical value $c$ of $F$, the permutation in $\Mon(P)$ corresponding to $c$ is a transposition.  
Since a primitive permutation group containing 
a transposition is a full symmetric group  (see \cite{wi}, Theorem 13.3), we conclude that \linebreak $\Mon(F)=S_m.$ \qed 

Let $F$ and $H$ be rational functions of degrees $n$ and $m$, and $H_1,$ $H_2$ and 
$F_1$, $F_2$ pairs of polynomials without common roots such that 
$H=H_1/H_2$ and $F=F_1/F_2.$ Let us define algebraic curves $h_{F,H}(x,y)$ and $h_F(x,y)$ by the formulas 
\be \l{cur} h_{H,F}:\ H_1(x)F_2(y)-H_2(x)F_1(y)=0,\ee 
and 
\be \l{cur2} h_F:\ \frac{F_1(x)F_2(y)-F_2(x)F_1(y)}{x-y}=0.\ee 
In case these curves are irreducible, their  genera can be calculated explicitly in terms of ramification of $H$ and $F$ as follows. 
Let $S=\{z_1, z_2, \dots , z_r\}$ be the union of all critical values of $H$ and $F$. 
For $i,$ $1\leq i \leq r,$ we denote by 
$$(a_{i,1},a_{i,2}, ... , a_{i,p_{i}})
$$ 
the collection of multiplicities of $H$ at the points of $H^{-1}\{z_i\}$, and by 
$$
(b_{i,1},b_{i,2}, ... , b_{i,q_i})
$$
the collection of multiplicities of $F$ at the points of $F^{-1}\{z_i\}$. 
 In this notation, 
the following formulas hold (see \cite{f3} or \cite{pq}): 
\be \l{rh0} 2-2g(h_{H,F})=\sum_{i=1}^{r}
\sum_{j_2=1}^{q_{i}}\sum_{j_1=1}^{p_{i}}\GCD(a_{i,j_1}b_{i,j_2})-mn(r-2),\ee
\be \l{rh1} 4-2g(h_{F})=\sum_{i=1}^{r}
\sum_{j_2=1}^{p_{i}}\sum_{j_1=1}^{p_{i}} \GCD(b_{i,j_1}b_{i,j_2})-(r-2)m^2.\ee

\bt \l{goo} 
Let $F$ be a simple rational function of degree $m\geq 4$,  and $H$ a  rational function of degree $n\geq 2$ such that  the curve $h_{H,F}$ is irreducible. Then 
$g(h_{H,F})>0$. In particular, the functional equation $F\circ X=H\circ Y$ has no solutions in rational functions $X,Y$. 
\et 
\pr Keeping the above notation, let us observe that if $z_i$, $1\leq i \leq r$, is not a critical value of $F$, then obviously 
\be \l{sth} \sum_{j_1=1}^{p_{i}}\GCD(a_{i,j_1}b_{i,j_2})=p_i, \ \ \  1\leq j_2\leq q_i,\ee
and
\be \l{urr1} \sum_{j_2=1}^{q_{i}}\sum_{j_1=1}^{p_{i}}\GCD(a_{i,j_1}b_{i,j_2})=mp_i.\ee 

Assume now that $z_i$, $1\leq i \leq r$, is a critical value of $F$. Then \eqref{sth} still holds  
if $b_{i,j_2}=1$, while if $b_{i,j_2}=2$, $1\leq j_2\leq q_i,$ we have  
$$\sum_{j_1=1}^{p_{i}}\GCD(a_{i,j_1}b_{i,j_2})=p_i+l_i,$$ where $l_i$ is the number of even numbers among the 
numbers $a_{i,j_1},$ $1\leq j_1\leq p_i.$ Since among the numbers $b_{i,j_2}$, $1\leq j_2\leq q_i,$  one number is equal to two and $m-2$ other numbers are equal to one, we conclude that
\be \l{urr2}  \sum_{j_2=1}^{q_{i}}\sum_{j_1=1}^{p_{i}}\GCD(a_{i,j_1}b_{i,j_2})=(m-2)p_i+p_i+l_i=mp_i+(l_i-p_i).\ee

As 
$$2n-2=\sum_{z\in \C\P^1}(\deg_zH-1)=\sum_{i=1}^{r}\sum_{j_1=1}^{p_{i}}(a_{i,j_1}-1)=rn-\sum_{i=1}^{r}{p_i},
$$
the equality 
\be \l{zxa} 
\sum_{i=1}^{r}{p_i}=(r-2)n+2
\ee
holds, implying by \eqref{urr1} and \eqref{urr2} that 
\be \l{itff}\sum_{i=1}^{r}
\sum_{j_2=1}^{q_{i}}\sum_{j_1=1}^{p_{i}}\GCD(a_{i,j_1}b_{i,j_2})=\sum_{i=1}^{r}mp_i+\sum{}^{'}(l_i-p_i)=\ee  
$$=m((r-2)n+2)+\sum{}^{'}(l_i-p_i),
$$ where the sum $\sum{}^{'}$ runs only over indices corresponding to critical values of $F$. 
It follows now from \eqref{rh0} that $g(h_{H,F})=0$ if and only if 
$$2m-2+\sum{}^{'}(l_i-p_i)=0.$$
Stating differently, $g(h_{H,F})=0$ if and only if the preimage 
$H^{-1}\{c_1,c_2,\dots ,c_{2m-2}\},$ where $c_1,c_2,\dots ,c_{2m-2}$ are critical values of $F$, contains exactly $2m-2$ points where the  multiplicity of $H$ is odd.

Let us observe now that for any finite subset $S$ of $\C\P^1$ it follows from  
$$2n-2=\sum_{z\in \C\P^1}(\deg_zH-1)\geq \sum_{z\in H^{-1}(S)}(\deg_zH-1)$$ that  the preimage $H^{-1}(S)$ contains at least $n(\vert S\vert -2)+2$ points and the equality is attained if and only if $S$ contains the set of critical values of $H$. 
Therefore, 
\be \l{pre} H^{-1}\{c_1,c_2,\dots ,c_{2m-2}\}\geq (2m-4)n+2,\ee 
and the equality is attained if and only if any critical value of $H$ is a critical value of $F$. On the other hand, the condition that $H^{-1}\{c_1,c_2,\dots ,c_{2m-2}\}$ contains $2m-2$ points where the  multiplicity of $H$ is odd implies that  
\be \l{pre2} H^{-1}\{c_1,c_2,\dots ,c_{2m-2}\}\leq (2m-2)+\frac{n(2m-2)-(2m-2)}{2}=(n+1)(m-1),\ee and the equality is attained if and only if all $2m-2$ points  with  odd
 multiplicity in 
 $H^{-1}\{c_1,c_2,\dots ,c_{2m-2}\}$ have multiplicity one, 
while all points with even multiplicity  have multiplicity two.    
Thus, if $g(h_{H,F})=0$, then 
$$(2m-4)n+2\leq (n+1)(m-1),$$ implying that $$(n-1)(m-3)\leq 0.$$
Since the last inequality is satisfied only for $n=1$ or for $m=2,3$, we conclude that $g(h_{H,F})>0$. \qed

\bt \l{goo1} 
Let $F$ be a simple rational function of degree $m\geq 3$. Then the curve $h_{F}$ is irreducible and  
$g(h_{F})>0$. In particular, the equality $F\circ X=F\circ Y$, where $X$ and $Y$ are rational functions, implies that $X=Y$. 
\et 
\pr 
It is well-known (see e.g. \cite{pq}, Corollary 2.3) that the curve $h_{F}(x,y)$ is irreducible if and only if the monodromy group $\Mon(F)$ is 
doubly transitive. Therefore,  since a symmetric group is doubly transitive, the irreducibility of $h_{F}(x,y)$ 
follows from Theorem \ref{0}.

Further, applying  \eqref{itff} and \eqref{zxa}  for 
 $H=F$  we see that 
 $$\sum_{i=1}^{r}
\sum_{j_2=1}^{p_{i}}\sum_{j_1=1}^{p_{i}}\GCD(b_{i,j_1}b_{i,j_2})=\sum_{i=1}^{2m-2}mp_i+\sum_{i=1}^{2m-2}(1-p_i)=\sum_{i=1}^{2m-2}(m-1)p_i+2m-2=$$
 $$=(m-1)((2m-4)m+2)+2m-2=m^2(2m-4)-2m^2+8m-4.$$ 
Therefore, by formula \eqref{rh1}, we have:  \be \l{burry} g(h_{F})=(m-2)^2,\ee and hence $g(h_{F})>0$ whenever $m\geq 3$. \qed

\subsection{Conditions for reducibility of $H(x)-F(y)=0$}
The  problem of finding conditions under which the algebraic curve $h_{H,F}$ is reducible, the so-called Davenport-Lewis-Schinzel problem, has a long story and  is not solved yet in its full generality (see \cite{fc} for an introduction to the topic). In this section, we consider a very particular case of this problem (Theorem \ref{go} below), which is related to the subject of this paper and can be handled without using serious group theoretic methods. The reader interested in these methods is referred to the recent paper \cite{nef2},  where a significant progress has been made in the polynomial case, and the bibliography therein.

Let $F$ be a rational function of degree $m\geq 2$, and  
$U\subset \C\P^1$ a simply connected domain containing no critical values of $F$. Then in $U$ there exist $m=\deg F$ different branches of the 
algebraic function $F^{-1}(z)$. We will denote these branches by small letters
$f_1,f_2,\dots ,f_m.$

\bl  \l{l1} Let $F$ be a  rational function of degree $m\geq 2$ such that $\Mon(F)=S_m,$ \linebreak 
 $f_1,f_2,\dots ,f_m$ different branches of $F^{-1}(z)$ defined in some simply connected domain $U$ containing no critical values of $F$, and  $C_i,$ $0\leq i \leq m,$ rational functions. 
Then the equality 
\be \l{eq} C_1 f_1 +C_2 f_2 +\dots +C_m f_m = C_0 \ee
 implies that $C_1=C_2= \dots = C_m.$
\el
\pr Assume, say, that $C_1\neq C_2$. Since $\Mon(F)=S_m$, the transposition $\sigma=(1,2)$ is contained 
in $\Mon(F)$, and considering the analytical continuation of equality \eqref{eq} 
along a loop corresponding to $\sigma$ we obtain the 
equality  \be \l{eq1} C_1  f_2  +C_2  f_1  +\dots +C_m  f_m  = C_0  .
\ee
It follows now from \eqref{eq} and \eqref{eq1}  that 
$$(C_1-C_2)(f_1-f_2)=0,$$ whence 
$f_1=f_2$ in contradiction with the assumption that $f_1,f_2,\dots ,f_m$ are diffe\-rent. \qed

\bl  \l{l2} Let $H$ be an indecomposable rational function of degree $n\geq 2$, 
$h_1,h_2,\dots ,h_m$ different branches of $H^{-1}(z)$ defined in some simply connected domain $V$ containing no critical values of $H$, and $R$ another rational function.  Then either  
$R(h_i)\neq R(h_j)$ for $i\neq j,$ $1\leq i, j \leq n$,  or $R(h_1)=R(h_2)=\dots =R(h_n)$. In the last case, $R(h_1)$ is a rational function.
\el
\pr 
It is easy to see that for fixed $j,$ $1\leq j \leq n,$ the set of all $i$, $1\leq i \leq n,$ such that 
$R(h_i)=R(h_j)$ is a block of $\Mon(H)$. Since $H$ is indecomposable, this implies the first  statement of the lemma.  Finally, if the functions $R(h_i),$ $1\leq i \leq n$, are equal, then the algebraic function obtained by a full analytical continuation of  
$R(h_1)$ is single-valued in $\C\P^1$ and therefore it is a rational function. \qed

\bt \l{go}  Let $H$ and $F$ be rational functions of degrees  
$n\geq 2$ and $m\geq 2$ such that $H$ 
is indecomposable, $\Mon(F)=S_m$,  and the curve $h_{H,F}$ is reducible. 
Then either $H=F\circ \mu,$ where $\mu$ is a M\"obius transformation, or 
$n=\binom{m}{k}$ for some $k,$ $1<k <m-1.$

\et

\pr Suppose that 
\be \l{fro} h_{H,F}(x,y)=M(x,y)N(x,y)\ee for some non-constant polynomials $M(x,y),$ $N(x,y).$ Notice that since  
$H$ and $F$ are non-constant, for such polynomials the degrees $\deg_xM,$ $\deg_yM,$ $\deg_xN,$ 
$\deg_yN$ are  distinct from zero.  

Let $h_1, h_2,\dots h_n$ be different branches of $H^{-1}(z)$ and $f_1, f_2,\dots f_m$ different branches of $F^{-1}(z)$ defined in a simply connected domain $U\subset \C\P^1$ containing no critical values of $F$ or $H.$ Since 
$$h_{H,F}(h_1,f_i)=0, \ \ \ 1\leq i \leq m,$$ among the indices $i,$ $1\leq i \leq m,$ there are $k=\deg_yM>0$ indicies for which 
the equality \be \l{eq2} M(h_1,f_i)=0\ee holds.  Moreover, $k<m$, since otherwise the equality $\deg_yM+\deg_yN=m$ implies that $\deg_yN=0.$ 
 Let $i_1,i_2,\dots ,i_k$ be indices for which \eqref{eq2} holds. Writing $M(x,y)$ in the form 
$$M(x,y)=P_k(x)y^k+P_{k-1}(x)y^{k-1}+\dots +P_1(x)y+P_0(x),$$ where $P_i,$ $0\leq i \leq k,$ are polynomials, we see that  \eqref{eq2} implies that 
\be \l{eq3} f_{i_1}+f_{i_2}+\dots + f_{i_k}=Q(h_1),\ee where $Q=P_{k-1}/P_k$ is a rational function. Furthermore, since the set $\{i_1,i_2,\dots ,i_k\}$ is a proper subset of  $\{1,2, \dots, m\}$, the function $Q(h_1)$ is not a rational function by Lemma \ref{l1}.  Therefore, by Lemma \ref{l2},  the functions $Q(h_i)$, $1\leq i \leq n,$ are pairwise different. 

Continuing equality \eqref{eq3} analytically  along an arbitrary  closed  curve $\gamma$ in $\C\P^1$, we obtain an equality where on the left side is a sum of branches of $F^{-1}(z)$ over a subset of $\{1,2, \dots, m\}$ containing $k$ elements, while on the right side is a branch $Q(h_i)$, $1\leq i \leq n$, of the function $Q(h_1)$.  
Furthermore, to different subsets of $\{1,2, \dots, m\}$  correspond different branches  of  $Q(h_1)$, for otherwise 
 subtracting we obtain a contradiction with Lemma \ref{l1}. Since 
the equality  $\Mon(F)=S_m$ implies that for  an appropriately chosen $\gamma$ we can obtain on the left side  a sum of branches of $F^{-1}(z)$ over any $k$-element subset of $\{1,2, \dots, m\}$, while the transitivity of 
 $\Mon(H)$ implies that for an appropriately chosen $\gamma$ we can obtain on the right side any branch  of $Q(h_i)$, $1\leq i \leq n,$ 
we conclude that $n$ is equal to the number of $k$-element subsets of $\{1,2, \dots, m\}$, for some $k,$ $1\leq k \leq n,$  that is, 
 $n=\binom{m}{k}$.  

To finish the proof, let us observe that if $k=1$, then $n=\binom{m}{k}$ implies that 
$n=m$. Furthermore, equality \eqref{eq3} implies the equality 
$$z=F(f_{i_1})=(F\circ Q)(h_1).$$ Thus, the function $F\circ Q$ is inverse to $h_1$, that is,  
   $H=F\circ Q$. Finally, $Q$ is a M\"obius transformation since $n=m$. The same conclusion  is true if $k=m-1$, since we can switch between $M(x,y)$ and $N(x,y)$. \qed

\subsection{Prime divisors of $\binom{m}{k}$} 
The classical theorem of 
Sylvester \cite{syl} and Schur \cite{sch} states that in the set of integers $a, a+1, \dots, a+b-1,$ where $a>b,$ there is a number divisible by a prime greater than $b$. 
For a natural number $x$, let us denote by $\f P(x)$ the greatest prime factor of $x$. Then  the  theorem of Sylvester and Schur may be reformulated as follows (\cite{erd}): for any  $m\geq 2k$ the inequality $\f P(\binom{m}{k})>k$ holds. Furthermore, the last inequality may be sharpened to the inequality \be \l{eh} \f P\left(\binom{m}{k}\right)\geq \frac{7}{5} k \ee (see \cite{fau}, \cite{han}). 
We will prove that this implies the following corollary.

\bt \l{gop} Let $m\geq 4$ be a natural number, and $k$ a natural number such that $1<k<m-1.$
Then there exists a prime number $p$ such that $p\mid \binom{m}{k}$
but $p\nmid m.$
\et

\pr 
Since $\binom{m}{k}=\binom{m}{m-k}$, it is enough to prove the theorem under the assumption that $m\geq 2k.$ Applying the Sylvester-Schur theorem, we conclude that there is a number 
$s$, $m-k+1\leq s \leq m$, such that $\f P(\binom{m}{k})=\f P(s)=p>k$. Moreover, if $s$ is  strictly less than $m$, then $p$ cannot be a divisor of $m$ for otherwise \be \l{f} p\vert (m-s),\ee where $m-s\leq k-1$ in contradiction with \be \l{o} p>k.\ee Since however $s$  can be {\it equal} to $m$, we modify slightly this argument. Namely, we apply the Sylvester-Schur theorem in its strong form \eqref{eh}  to the binomial coefficient $\binom{m-1}{k-1}$ related with $\binom{m}{k}$ by the equality 
\be \l{s} \binom{m}{k}=\frac{m(m-1)\dots (m-k+1)}{k(k-1)\dots 1}=\frac{m}{k}\binom{m-1}{k-1}.\ee
Notice that \eqref{s} implies that every prime factor $p$ of $\binom{m-1}{k-1}$
satisfying  \eqref{o}  remains a prime factor of $\binom{m}{k}$.

 Since
$m\geq 2k$ implies $m-1\geq 2(k-1)$, applying \eqref{eh} to $\binom{m-1}{k-1}$ we conclude that there is a number 
$s$, $m-k+1\leq s \leq m-1$, such that
$$\f P\left(\binom{m-1}{k-1}\right)= \f P(s)=p\geq\frac{7}{5}(k-1).$$ Furthermore, if $k>3$ then $$p\geq \frac{7}{5}(k-1)>k,$$ implying that  $p\mid \binom{m}{k}.$
On the other hand, $p\nmid m$ since otherwise  \eqref{f} holds in contradiction with \eqref{o}. 

For $k\leq 3$, the theorem can be proved by an elementary argument. If $k=2$, then $$\binom{m}{k}=\frac{m(m-1)}{2}.$$ Therefore, since $\GCD(m,m-1)=1$, the statement of the theorem is true, whenever  
$(m-1)\nmid 2$, and the last condition is always satisfied if $m>3.$  
Similarly, if $k=3$, then $$\binom{m}{k}=\frac{m(m-1)(m-2)}{2 \cdot 3},$$  and the statement of the theorem is true whenever $(m-1)\nmid 6.$ The last condition  fails to be true for $m>3$ only if $m$ is equal to  $4$ or $7$. However, the pair $m=4$, $k=3$ does not satisfy the condition $1<k< m-1.$ On the other hand, for the pair $m=7$, $k=3$, we have $\binom{m}{k}=\binom{7}{3}=5\cdot 7$, and the statement of the theorem is satisfied for $p=5.$ \qed

\vskip 0.2cm
\noindent{\it Proof of Theorem \ref{t1}.} The proof is by induction on $l.$ 
Let 
\be \l{en} F^{\circ l} =F_r\circ F_{r-1}\circ \dots \circ F_1\ee
be a decomposition of $F^{\circ l}$, $l\geq 1$, into a composition of indecomposable  rational functions. Since $F$ is indecomposable by Theorem \ref{0}, for $l=1$ the theorem  is true. On the other hand, since by Theorem \ref{goo1} the equality $F\circ X=F\circ Y$ implies that $X=Y$, to prove the inductive step it is enough to show that    
equality \eqref{en} implies that 
\be \l{mu} F_r=F\circ \mu \ee for some M\"obius transformation $\mu.$

Clearly, equality \eqref{en} implies that the algebraic curve 
\be \l{us} F(x)-F_r(y)=0\ee has a factor of genus zero. Therefore, \eqref{us} is  reducible by Theo\-rem \ref{goo}. Since $\Mon(F) =S_m$ by Theorem \ref{0},  it follows now from  Theorem \ref{go} that either \eqref{mu} holds, or 
$\deg F_r=\binom{m}{k}$ for some $k,$ $1<k <m-1.$ However, the last case is impossible since 
 \eqref{en} implies that any prime divisor of $\deg F_r$ is a prime divisor of $\deg F$ in contradiction with Theorem \ref{gop}. 
  \qed

\bc \l{cor1} Let $F$ be a simple rational function of degree $m\geq 4$, and $G_i$, $1\leq i \leq r,$  rational functions of degree at least two such that  $$F^{\circ l}=G_r\circ G_{r-1}\circ \dots G_1$$ for some $l\geq 1$. Then there exist M\"obius transformations $\nu_i$, $1\leq i <r,$ and integers $s_i\geq 1,$ $1\leq i \leq r,$
such that 
$$G_r=F^{\circ s_r}\circ \nu_{r-1}, \ \ \ 
G_i=\nu_{i}^{-1}\circ F^{\circ s_i} \circ \nu_{i-1}, \ \ \ 1<i< r, \ \ \ {\it and} \ \ \ G_1=\nu_{1}^{-1}\circ F^{\circ s_1}.$$ 
\ec 
\pr To prove the corollary, it is enough to  
decompose each $G_i$, $1\leq i \leq r,$ into a composition of indecomposable rational functions and to apply 
Theorem \ref{t1}. \qed

\section{Groups and semigroups related to rational functions}
\subsection{\l{rpr} Groups and semigroups related to simple rational functions}
 We start this section by recalling some basic facts concerning the groups and semigroups 
defined in the introduction. We will say that a  rational function $F$ of degree at least two is {\it special} if  $F$ is either a Latt\`es map, or it is conjugate to a power $z^{\pm n}$  or to a Chebyshev polynomial $\pm T_n$  (we will recall the definition of  Latt\`es maps below, in Section \ref{44}, in a more general context). 

It is obvious that $C(F)$ is a semigroup, and it follows from the inclusions  
\be \l{lcm} C(F^{\circ k}),\  C(F^{\circ l})\subseteq C(F^{\circ \LCM(k,l)})\ee 
that $C_{\infty}(F)$ is also a semigroup.
We  use the following characterization of $C_{\infty}(F)$. 

\bt \l{rr} Let $F$ be a non-special rational function  
of degree at least two. Then a rational function    $G$ 
of degree at least two belongs to $C_{\infty}(F)$ if and only if equality \eqref{u} holds for some $k,l\geq 1$.
\et 
\pr 
By the Ritt theorem (see \cite{r}, and also \cite{e2}, \cite{revi}), commuting rational functions  
of degree at least two are either special or have a common iterate. Thus, if $G$ belongs to $C_{\infty}(F)\setminus \Aut_{\infty}(F)$, then  \eqref{u} holds for some $k,l\geq 1$. On the other hand, if \eqref{u} holds, then $G$ commutes with 
$F^{\circ k}$, and thus $G$ belongs to $C_{\infty}(F).$ 
\qed 
 
Notice that a practical method for describing $C(F)$ for an arbitrary non-special rational function $F$ was given in the recent paper \cite{revi}, but  
a satisfactory description of $C_{\infty}(F)$ 
is still not known (see \cite{p20} for some particular results). Thus, condition \eqref{u} remains  the only characterization of $C_{\infty}(F)$.

Let us recall that by the results of Freire, Lopes, Ma\~n\'e (\cite{flm}) and Lyubich (\cite{l}), 
 for every rational function $F$ of degree $n\geq 2$, there exists a unique probability measure $\mu_F$ on $\C\P^1$, which is invariant under $F$, has support equal to the Julia set $J_F$, and achieves maximal entropy $\log n$ among all $F$-invariant probability measures.  The measure $\mu_F$ can be characterized as follows. For $a\in \C\P^1$,  
let  $z_i^k(a),$ \linebreak $i=1, \dots, n^k,$ be the roots of the equation
$F^{\circ k}(z) = a$ counted with multiplicity, and $\mu_{F,k}(a)$ be  
the measure defined by 
\be \l{mera} \mu_{F,k}(a)=\frac{1}{n^k}\sum_{i=1}^{n^k}\delta_{z_i^k(a)}.\ee Then 
for every  $a\in \C\P^1$ with two possible exceptions, the sequence  $\mu_{F,k}(a)$, $k\geq 1$,  converges in the weak topology to $\mu_F.$ It follows from the characterization of $\mu_F$ as a limit of \eqref{mera}  that any $G$ sharing an iterate with $F$ belongs to $E(F)$. Thus,  
\be \l{new} C_{\infty}(F)\subseteq E(F).\ee Moreover, since the equality 
$$F^{\circ n}=\alpha^{-1} \circ F^{\circ n}\circ \alpha,$$ where $\alpha\in \Aut(\C\P^1)$ and $n\geq 1$, implies that 
$$\vert S\cap F^{-nk}(a)\vert =\vert \alpha(S)\cap F^{-nk}(\alpha(a))\vert, \ \ \ k\geq 1, $$ for any set $S\subset \C\P^1$ and $a\in \C\P^1$, 
this characterization   yields that 
\be \l{yby} \Aut_{\infty}(F)\subseteq E_0(F).\ee  

The fact that $E(F)$ is a semigroup can be established using the Lyubich operator or the balancedness property of $\mu_F$ (see \cite{peter}, \cite{p20}), and the analogue of  Theorem \ref{rr} is the following statement.

\bt \l{rrr} Let $F$ be a non-special rational function  
of degree at least two. Then a rational function    $G$ 
of degree at least two belongs to $E(F)$ if and only if equality \eqref{uu} holds for some $k_1,l\geq 1$, $k_2\geq 0.$
\et 
\pr It is known and can be easily shown using the Lyubich operator that equality \eqref{uu} implies the equality $\mu_F=\mu_G$ (\cite{lev}). On the other hand, it is shown in \cite{ye} that the characterization of non-special rational functions sharing the measure of maximal entropy obtained in the papers \cite{lev}, \cite{lp} implies that for such functions $F$ and $G$ equality 
\eqref{uu} holds. \qed 

A complete description of $E(F)$ is known 
only if $F$ is a polynomial, in which case $E(F)\setminus E_0(F)$ coincides with the set of polynomials sharing a Julia set with $F$ (see \cite{a2}, \cite{b}, \cite{sh} and also \cite{sm}, \cite{p20}). Some partial results in the rational case can be found in  \cite{entr}, \cite{ye}.

The group $G_0(F)$ obviously is a subgroup of the larger group 
$ G(F)$ consisting of all $\sigma\in\rm{Aut}(\C\P^1)$ such that \be \l{eblys} F\circ \sigma=\nu\circ F\ee for some  $\nu\in\rm{Aut}(\C\P^1)$. 
 It is easy to see that 
 $ G(F)$ is indeed a group 
and that the map 
\be \l{homic}\gamma:\sigma \rightarrow \nu_{\sigma}\ee is a homomorphism from   $ G(F)$ to the  group 
$\rm{Aut}(\C\P^1)$. The group $G(F)$ is finite and its order is  bounded in terms of $m=\deg F$,  unless \be \l{for} \alpha\circ F\circ \beta=z^m\ee for some M\"obius transformations $\alpha,\beta$ (see \cite{fin}, Section 4 or \cite{sym}, Section 2). 
Thus, the group $G_0(F)$ is also finite, unless \eqref{for} holds.

\bl \l{g} 
Let $F$ be a simple rational function of degree $m\geq 3$. Then  the group  $G_0(F)$ is finite, and the restriction  of $\gamma$ to $G_0(F)$  is an automorphism of $G_0(F)$.
\el
\pr Since equality \eqref{for} is impossible for simple $F$ of degree $m\geq 3$,  the group  $G_0(F)$ is finite. Furthermore, since the equality $F=F\circ \sigma$, where $\sigma\in  G_0(F)$, implies by Theo\-rem \ref{goo1} that $\sigma$ is the identity element,  
the restriction of $\gamma$ on $G_0(F)$ is a monomorphism. Since $\gamma(G_0(F))\subseteq G_0(F)$ by the definition, this implies that  the restriction  of $\gamma$ to $G_0(F)$ 
is an automorphism of $G_0(F)$. \qed 

\bc \l{svin1} Let $F$ be a simple rational function of degree $m\geq 3$. Then \linebreak $G_0(F)\subseteq \Aut(F^{\circ s}),$ where $s=\vert \Aut(G_0(F))\vert .$  
\ec
\pr 
For  $s=\vert \Aut(G_0(F))\vert $, the iterate $\gamma^{\circ s}$ is the identity automorphism of $G_0(F)$. Therefore, since for every $\sigma\in G_0(F)$ the equality 
$$F^{\circ l}\circ \sigma=\gamma^{\circ l}(\sigma)\circ F^{\circ l},\ \ \  l\geq 1,$$  holds, 
 every $\sigma\in G_0(F)$ commutes with $F^{\circ s}$.  \qed

\bl \l{73} Let $F$ be a rational function, and $\sigma$ a M\"obius transformation such that 
 \be \l{xer} (\sigma \circ F )^{\circ l}=F^{\circ l}\ee 
for some $l\geq 1.$ Then $\sigma\in \Aut(F^{\circ l}).$ 

\el  
\pr Clearly, equality \eqref{xer} 
 implies   the equality \be \l{egas} (\sigma \circ F )^{\circ (l-1)}\circ \sigma=F^{\circ (l-1)}. \ee
Composing now $F $ with the both parts of equality \eqref{egas}, we obtain the equality 
\be \l{xur} ( F \circ\sigma )^{\circ l}=F^{\circ l}.\ee 
It follows now from \eqref{xer} and \eqref{xur} that  
$$F^{\circ l}\circ \sigma=(\sigma \circ F )^{\circ l}\circ \sigma= \sigma\circ ( F \circ\sigma )^{\circ l}=\sigma\circ F^{\circ l}. \eqno{\Box}$$

\bl \l{lats0} Let $F$ be a simple rational function of degree $m\geq 4$. Then $F$ is not a special function. 
\el 
\pr The proof follows easily from the analysis of  ramifications of special functions. Since below we prove a more general result (Lemma \ref{lats}), we omit it.  
\qed

\bt \l{32} Let $F$ be a simple rational function of degree $m\geq 4$. Then 
$$E_0(F)=G_0(F)=\Aut_{\infty}(F)=\Aut(F^{\circ s}),$$ where $s=\vert \Aut(G_0(F))\vert .$ 
\et 
\pr By  Corollary \ref{svin1} and \eqref{yby}, we have:
$$G_0(F)\subseteq \Aut(F^{\circ s})\subseteq \Aut_{\infty}(F)\subseteq E_0(F).$$
Thus, to prove the theorem we only must  prove that $E_0(F)
\subseteq G_0(F).$ For this, it is enough to show that 
for every  $\sigma \in E_0(F)$ there exists $\nu \in E_0(F)$ such that \eqref{eblys} holds. Let $\sigma$ be an element of $E_0(F).$ Then  $F\circ \sigma$ is a simple rational function which belongs to $ E(F)$, implying by Lemma \ref{lats0} and Theorem \ref{rrr} that 
$$ F^{\circ k_1}=F^{\circ k_2}\circ (F\circ \sigma )^{\circ l}$$   for some $k_1,l\geq 1$, $k_2\geq 0.$ Applying to the last equality recursively Theorem \ref{goo1}, we see that
\be \l{sinn}  F^{\circ (k_1-k_2)}=(F\circ \sigma )^{\circ l}.\ee
Therefore, by  Theorem \ref{t1}, 
there exist M\"obius transformations  $\mu_i,$ $1\leq i \leq l-1,$ such that 
\be \l{eba} F\circ \sigma =F \circ \mu_{r-1}, \ \ \ 
F\circ \sigma =\mu_{i}^{-1}\circ F  \circ \mu_{i-1}, \ \ \ 1<i< l, \ \ \ {\rm and} \ \ \ F\circ \sigma =\mu_{1}^{-1}\circ F. 
\ee  Thus, equality \eqref{eblys} holds 
 for $\nu=\mu_{1}^{-1}.$ 
Furthermore, since equality \eqref{sinn} implies that $k_1-k_2=l$, we have:
$$F^{\circ l}=(F\circ \sigma)^{\circ l}=(\mu_{1}^{-1}\circ F)^{\circ l},$$ 
implying by Lemma \ref{73} that $$\mu_{1}^{-1}\in \Aut(F^{\circ l})\subseteq \Aut_{\infty}(F)\subseteq E_0(F).\eqno{\Box}$$

\bc \l{svin2} Let $F$ be a simple rational function of degree $m\geq 4$.   
Then  every element of the semigroup $\langle \Aut_{\infty}(F), F \rangle$ can be represented in a unique way in the form $\alpha\circ F^{\circ s}$, where  $\alpha\in \Aut_{\infty}(F)$ and $s\geq 0$. Moreover, for every  $k\geq 1$, every element of the semigroup $\langle \Aut(F^{\circ k}), F \rangle$ can be represented in a unique way in the form $\alpha\circ F^{\circ s}$, where  $\alpha\in \Aut(F^{\circ k})$ and $s\geq 0$. 
\ec
\pr The first part of the corollary follows from the equality $\Aut_{\infty}(F)=G_0(F).$ 
 To prove the second, it is enough to observe that 
for every $\nu \in \Aut(F^{\circ k})$ the element  
  $\nu'\in \Aut_{\infty}(F)$  such that \be \l{intu} F\circ \nu=\nu'\circ F\ee belongs to $ \Aut(F^{\circ k})$.   Indeed, \eqref{intu} implies that  
$$F^{\circ k}\circ \nu'\circ F=F^{\circ k}\circ F \circ \nu=F\circ F^{\circ k} \circ \nu=F\circ  \nu\circ F^{\circ k} =\nu'\circ F\circ F^{\circ k}=\nu'\circ F^{\circ k}\circ F,$$ whence  $\nu' \in \Aut(F^{\circ k})$. \qed

\vskip 0.2cm

\noindent{\it Proof of Theorem \ref{t2}.} 
In view of Theorem \ref{32}, we only must show that  
$$C_{\infty}(F)=E(F)=\langle \Aut_{\infty}(F), F \rangle. $$
By \eqref{new},  the first equality follows from Theorem \ref{rrr} and  Theorem \ref{goo1}, since the latter implies that any $G$ satisfying \eqref{uu} satisfies  \eqref{u} for $k=k_1-k_2$. Since the semigroup $\langle \Aut_{\infty}(F), F \rangle$ is obviously a subsemigroup of $C_{\infty}(F)$, to finish the proof we only must show that if a rational function $G$ satisfies \eqref{u}, then it belongs to $\langle \Aut_{\infty}(F), F \rangle$.

 Applying  Corollary \ref{cor1} to equality \eqref{u}, we see that 
there exist M\"obius transformations  $\mu_i,$ $1\leq i \leq l-1,$ such that 
$$ G=F^{\circ s}\circ \mu_{l-1}, \ \ \ 
G=\mu_{i}^{-1}\circ F^{\circ s} \circ \mu_{i-1}, \ \ \ 1<i< l, \ \ \ {\rm and} \ \ \ G=\mu_{1}^{-1}\circ F^{\circ s},
$$ where $s=k/l$. 
 Moreover, since 
$$F^{\circ sl}=G^{\circ l}=(\mu_{1}^{-1}\circ F^{\circ s})^{\circ l},$$ 
Lemma \ref{73} implies that $$\mu_{1}^{-1}\in \Aut(F^{\circ sl})\subseteq \Aut_{\infty}(F).$$  Thus,   $$G=\mu_{1}^{-1}\circ F^{\circ s}\in \langle \Aut_{\infty}(F), F \rangle. \eqno{\Box} $$

\subsection{\l{gas} Groups and semigroups related to general rational functions} 
Let us recall that writing a rational function $F=F(z)$ of degree $m$ as 
$F=P/Q$, where
$$P(z)=a_mz^m+a_{m-1}z^{m-1}+...+a_1z+a_0, \ \ \ Q(z)=b_mz^m+b_{m-1}z^{m-1}+...+b_1z+b_0$$ 
are polynomials of degree $m$ without common roots, we can  
identify the space
of rational functions of degree $m$ with the algebraic variety $${\rm Rat}_m=\C\P^{2m+1}\setminus {\rm Res}_{m,m,z}(P,Q),$$ where ${\rm Res}_{m,m,z}(P,Q)$ denotes the resultant of $P$ and $Q$.  We recall that we say that some statement holds for general rational functions of degree $m$, if it holds for all \linebreak $F\in {\rm Rat}_m$ with exception of some proper Zariski closed subset.

\bl \l{lom1} A general rational function $F$ of degree $m\geq 2$ is simple. 
\el
\pr 
We recall that for $F\in {\rm Rat}_m$ the set of finite critical points of $F$ coincides with the set of zeroes of its Wronskian $$W(z)=P'(z)Q(z)-P(z)Q'(z).$$ Obviously, $\deg W\leq 2m-2$. Moreover, $\deg W=2m-2,$  unless $F$ belongs to the projective hypersurface $U$ in $\C\P^{2m+1}$ defined by 
$$U:\ a_mb_{m-1}-b_ma_{m-1}=0.$$ 

Let us define now a polynomial $R(t)$ by the formula 
$$R(t)={\rm Res}_{2m-2,m,z}(W(z),P(z)-Q(z)t).$$ By the well-known property of the resultant, for $F\in  {\rm Rat}_m\setminus U$ the equality 
\be \l{fo} R(t)=c\prod_{\zeta, W(\zeta)=0}(P(\zeta)-Q(\zeta)t)\ee holds for some $c\in \C^*$, and thus the set of zeroes of $R(t)$  coincides with  the  set of finite critical values of $F$. 

Finally, let us define a projective hypersurface $Z$ in $\C\P^{2m+1}$  by   
$$Z:\, {\rm Res}_{2m-2,2m-3,t}(R(t),R'(t))=0.$$ By the resultant properties,   
 $F\in {\rm Rat}_m\setminus U$ belongs to $Z$ if and only if either some finite critical values of $F$ coincide, or $\deg R(t)<2m-2$ meaning that infinity is a critical value of $F$. Thus,  
 every rational function $F\in {\rm Rat}_m\setminus Z\cup U$ has $2m-2$ distinct finite critical values, and hence  is simple.   \qed

\bl \l{lom2} For a general rational function $F$ of degree $m\geq 3$ the group $G(F)$ is trivial.  
\el
\pr Let $$\alpha=\frac{\alpha_{1,1}z+\alpha_{0,1}}{\alpha_{1,2}z+\alpha_{0,2}}, \ \ \ \beta=\frac{\beta_{1,1}z+\beta_{0,1}}{\beta_{1,2}z+\beta_{0,2}}$$ be elements of 
${\rm Rat}_2$, and 
$$F=\frac{f_{m,1}z^m+f_{m-1,1}z^{m-1}+\dots +f_{1,1}z+f_{0,1}}{f_{m,2}z^m+f_{m-1,2}z^{m-1}+\dots +f_{1,2}z+f_{0,2}}$$ an element of 
${\rm Rat}_m$. 
It is easy to see that the coefficients of the numerator of the rational function
$\alpha\circ F\circ \beta-F$ are polynomials in  $\alpha_i^j$, $\beta_i^j$,  $f_i^j$ homogenous of degree one in $\alpha_i^j$,  homogenous of degree two in $f_i^j$, and  homogenous of degree $m$ in $\beta_i^j$. Thus, the equality  
$$\alpha\circ F\circ \beta=F$$ implies that the coefficients of $\alpha,$ $\beta$, and $F$ belong to some projective algebraic variety $$W\subseteq \C\P^3\times \C\P^{2m+1} \times \C\P^3.$$ 
Since the projection $$p_k: \  (\C\P^1)^l\times (\C\P^1)^k\rightarrow 
(\C\P^1)^k, \ \ \ \ \ \ \ k,l\geq 1,$$ is a closed map (see e.g. \cite{mum}), this implies that the set of $F\in {\rm Rat}_m$ with non-trivial group $G(F)$   is contained in some Zariski closed subset of $\C\P^{2m+1}.$ 

To finish the proof, we only must show that $Z$ does not contain the whole variety ${\rm Rat}_m$ for $m\geq 3.$ For this, it is enough to show that for every $m\geq 3$ there exists a polynomial $F$ of degree $m$ such that the group $G(F)$ is trivial. Let us recall that for any polynomial $F$ of degree $m\geq 2$  the group 
$G(F)$ is a finite cyclic group generated by a polynomial, unless \eqref{for} holds (see e.g. \cite{sym},  Section 2). On the other hand, it is easy to see that if $F$ has the form \be \l{exx} F=z^m+a_{m-2}z^{m-2}+a_{m-3}z^{m-3}+\dots + a_0,\ee then equality \eqref{eblys} may hold for polynomials $\sigma=az+b$, $\mu=cz+d$ only if $b=0$ 
and $a$ is a root of unity. This implies easily that for any polynomial of the form \eqref{exx} 
with $a_{m-2}\neq 0,$ $a_{m-3}\neq 0$ the group $G(F)$ is trivial.   
\qed

Notice that Lemma \ref{lom2} is not true for $m=2.$ Indeed, for every rational function $F$ of degree two  there exist  M\"obius transformations
 $\alpha,\beta$ such that equality \eqref{for} holds, implying that  the group $G(P)$ is non-trivial,  and even infinite.

\vskip 0.2cm
\noindent{\it Proof of Theorem \ref{t3}.} Since $G_0(F)$ is a subgroup of $G(F),$ the theorem follows  from Theorem \ref{t2} combined with Lemma \ref{lom1} and Lemma \ref{lom2}. \qed

\section{\l{4} Semiconjugate rational functions and invariant curves} 
\subsection{\l{44} Generalized Latt\`es maps, semiconjugate rational functions, and invariant curves}
In this section, we recall some definitions and results related to  the functional equation 
\be \l{ax} A\circ X=X\circ B\ee in rational functions, 
and to  invariant curves for endomorphisms of $(\C\P^1)^2$ of the form 
 $$(A_1,A_2):\,  (z_1,z_2)\rightarrow (A_1(z_1),A_2(z_2)),$$
where  $A_1$, $A_2$ are rational functions.

Let us recall that an {\it orbifold} $\f O$ on $\C\P^1$  is a ramification function $\nu:\C\P^1\rightarrow \mathbb N$, which takes the value $\nu(z)=1$ except at a finite set of points. 
For an orbifold $\f O$, the set of {\it singular points} of $\f O$ is the set 
$$c(\f O)=\{z_1,z_2, \dots, z_s, \dots \}=\{z\in \C\P^1 \mid \nu(z)>1\},$$ and 
 the {\it  Euler characteristic} of $\f O$ is the number
$$ \chi(\f O)=2+\sum_{z\in R}\left(\frac{1}{\nu(z)}-1\right).$$

Let $A$ be a rational function, and $\f O_1$, $\f O_2$  orbifolds with ramifications functions $\nu_1$ and $\nu_2$. We say that  
$A:\,  \f O_1\rightarrow \f O_2$ is  {\it a covering map} 
between orbifolds
if for any $z\in \C\P^1$ the equality 
$$ \nu_{2}(A(z))=\nu_{1}(z)\deg_zA$$ holds.
In case  
the weaker condition 
$$ \nu_{2}(A(z))=\nu_{1}(z)\GCD(\deg_zA, \nu_{2}(A(z))$$ is satisfied,  
we say that $A:\,  \f O_1\rightarrow \f O_2$ is  {\it a  minimal holomorphic  map} 
between orbifolds. 

In the above terms, {\it a  Latt\`es map} can be defined as a rational function $A$  such that $A:\f O\rightarrow \f O$ is a  covering map
for some orbifold $\f O$ (see \cite{mil2}, \cite{lattes}). 
Following \cite{lattes}, we say 
 that a rational function $A$ of degree at least two is {\it a generalized Latt\`es map}  if there exists an orbifold $\f O$, 
 distinct from the non-ramified sphere, such that  $A:\,  \f O\rightarrow \f O$ is  a  minimal holomorphic  map. Thus, $A$ is a  Latt\`es map if there exists an orbifold $\f O$ such that 
\be \l{sin} \nu(A(z))=\nu(z)\deg_zA, \ \ \ \ z\in \C\P^1,\ee
and $A$ is  a  generalized Latt\`es map if there exists an orbifold $\f O$ such that 
\be \l{ee} \nu(A(z))=\nu(z)\GCD(\deg_zA,\nu (A(z))),  \ \ \ \ z\in \C\P^1.\ee 
Since \eqref{sin} implies \eqref{ee}, any Latt\`es map is a generalized Latt\`es map. 
More generally, any special function is a generalized Latt\`es map (see \cite{lattes}). 
 Notice  that if  \eqref{sin} holds for some rational function $A$ and orbifold
$\f O$, then $\f O$ has zero Euler characteristic, while \eqref{ee} implies that the Euler characteristic of $\f O$ is non-negative 
(see e.g. \cite{lattes} for more detail). 

Latt\`es maps and generalized Latt\`es maps can be characterized also in different terms (see \cite{mil2}, \cite{lattes}). However, the definition using orbifolds is most convenient for our purposes since it permits to show easily that a simple rational function of degree at least four is not a generalized Latt\`es map. In turn, this fact is crucial for our proof of Theorem \ref{t4} and Theorem  \ref{t5}  since for rational functions $A$ and $A_1,A_2$ that are not generalized Latt\`es maps describing solutions  of \eqref{ax} and $(A_1,A_2)$-invariant curves reduces to describing decompositions of iterates of $A$ and $A_1,A_2$.  

Specifically, our proof of Theorem \ref{t4} relies on the following corollary of the classification of semiconjugate rational functions (see  \cite{ic}, Proposition 3.3).

\bt \l{tt} Let $A,B$ be rational functions of degree at least two and $X$  a rational function of degree at least one such that equality \eqref{ax} holds and 
$A$ is not a generalized Latt\`es map. Then there exists a rational function $Y$ such that $X\circ Y=A^{\circ d}$ for some $d\geq 0.$
 \qed
\et

In turn, our proof of Theorem \ref{t5} uses the following corollary of the description of invariant curves 
for endomorphisms  $(A_1,A_2)$ of $(\C\P^1)^2$ 
(see \cite{ic}, Theorem 1.1).


\bt \l{1} Let $A_1$, $A_2$ be rational functions of degree at least two that are not generalized Latt\`es maps, and 
$C$ an irreducible algebraic  $(A_1,A_2)$-invariant  curve  in $(\C\P^1)^2$ that is not a vertical or horizontal line. Then
 there exist  rational functions $X_1,$ $X_2,$ $Y_1,$ $Y_2,$ $B$  such that: 
\begin{enumerate} 
\item[1.]
 The diagram 
\be  \l{xx}
\begin{CD} 
(\C\P^1)^2 @>(B,B)>>(\C\P^1)^2 \\ 
@V (X_1,X_2)  VV @VV  (X_1,X_2) V\\ 
 (\C\P^1)^2 @>(A_1,A_2)>> (\C\P^1)^2
\end{CD}
\ee
commutes, \item[2.]  The 
equalities 
\be \l{en1} X_1\circ Y_1=A_1^{\circ d}, \ \ \ \ \ \ \ 
 X_2\circ Y_2=A_2^{\circ d},\ee hold for some $d\geq 0$, \item[3.]  The map $t\rightarrow (X_1(t),X_2(t))$ is a parametrization of  $C$. \qed 
\end{enumerate}
\et

Notice that if diagram \eqref{xx} commutes, then this condition  alone obviously is sufficient for $(A_1,A_2)$-invariance of the curve $C$ parametrized by  $t\rightarrow (X_1(t),X_2(t))$.

\subsection{Proof of Theorem \ref{t4}, Theorem \ref{t5}, and Theorem \ref{t6}.}
We start by proving the following lemma. 

\bl \l{lats} Let $F$ be a simple rational function of degree $m\geq  4$. Then $F$ is not a generalized Latt\`es map. 
\el 
\pr If $F$ is a simple rational function of degree $m\geq  4$, then the preimage of any $k$ distinct points  of $\C\P^1$ under $F$   contains at least $k(m-2)\geq 2k$ distinct points $z$ such that $\deg_zF=1$. Thus, if the equality 
$$\nu(F(z))=\nu(z)\GCD(\deg_zF,\nu (F(z))),  \ \ \ \ z\in \C\P^1$$ 
 holds for some orbifold $\f O$ distinct from the non-ramified sphere, then the preimage $F^{-1}\{c(\f O)\}$ must contain at least  $2\vert c(\f O)\vert $ points where $\nu(z)>1$. However,
this is impossible since any such a point belongs to  $c(\f O)$. \qed

\vskip 0.2cm

\noindent {\it Proof of Theorem \ref{t4}.} By Lemma \ref{lats}, $F$ is not a generalized Latt\`es map. Since a rational function $F$ is a generalized Latt\`es map if and only if some iterate $F^{\circ d}$, $d\geq 1$, is a generalized Latt\`es map (see \cite{ic}, Section 2.3), this implies that $F^{\circ r}$ also is not a generalized Latt\`es map. Hence, by Theorem \ref{tt}, there exists a rational function $Y$ such that the equality  
$$X\circ Y=F^{\circ rd}$$ holds for some $d\geq 0$. By Corollary \ref{cor1}, this implies that 
$$X=F^{\circ l}\circ \mu$$ 
for some M\"obius transformation $\mu$ and $l\geq 0.$ Thus, diagram \eqref{ii2} reduces to the equality 
$$ F^{\circ r}\circ F^{\circ l}\circ \mu= F^{\circ l}\circ \mu\circ G,$$ 
and applying to this equality Theorem \ref{goo1}, we conclude that 
$$G=\mu^{-1}\circ F^{\circ r}\circ \mu.\eqno{\Box}$$

\vskip 0.2cm

\noindent{\it Proof of Theorem \ref{t5}.}  Assume that \be \l{lass} (F_1,F_2)^{\circ d}(C)=C, \ \ \ \ d\geq 1.\ee Then Theorem \ref{1} and Theorem \ref{t4} imply  that $C$ is parametrized by 
\be \l{buran} t\rightarrow \left((F_1^{\circ d_1}\circ \beta)(t),\,(F_2^{\circ d_2}\circ \alpha)(t)\right)\ee for some $d_1,d_2\geq 0$ and  M\"obius transformations $\alpha$, $\beta$  such that 
$$\beta^{-1}\circ F_1^{\circ d}\circ \beta=  \alpha^{-1}\circ F_2^{\circ d}\circ \alpha.$$  

 It is clear that without loss of generality we  may assume that $\beta$ is the identity map implying that 
\be \l{kosh1} F_1^{\circ d}=  \alpha^{-1}\circ F_2^{\circ d}\circ \alpha=(\alpha^{-1}\circ F_2\circ \alpha)^{\circ d}.\ee This yields that $$\alpha^{-1}\circ F_2\circ \alpha\in C(F_1^{\circ d})\subseteq  C_{\infty}(F_1),$$ and hence 
\be \l{kosh2} \alpha^{-1}\circ  F_2\circ  \alpha=\mu\circ F_1\ee for some $\mu\in \Aut_{\infty}(F_1)$ by Theorem \ref{t2} and Corollary \ref{svin2}. Further, equalities \eqref{kosh1} and \eqref{kosh2} imply by Lemma \ref{73} that $\mu\in \Aut(F_1^{\circ d})$. 
Therefore, 
\be \l{reby} F_2=\alpha\circ \mu\circ F_1\circ  \alpha^{-1}\ee for some $\mu\in \Aut(F_1^{\circ d}),$  and parametrization \eqref{buran} takes the form 
\be \l{krott} t\rightarrow \left(F_1^{\circ d_1}(t),\,\alpha\circ (\mu\circ F_1)^{\circ d_2}(t)\right).\ee
Moreover, it follows from \eqref{krott} by Corollary \ref{svin2} that there exists $\mu'\in  \Aut(F_1^{\circ d})$ such that this parametrization can be written in the form 
\be \l{krota} t\rightarrow \left(F_1^{\circ d_1}(t),\, (\alpha\circ\mu'\circ F_1^{\circ d_2})(t)\right).\ee 

If $d_1\leq d_2,$ then \eqref{krota} implies that $C$ is parametrized 
by  \be \l{sl1} t\rightarrow \left(t,(\alpha\circ \mu'\circ F_1^{\circ (d_2-d_1)})(t)\right)\ee for some  $\mu'\in \Aut(F_1^{\circ d}).$   
On the other hand, if $d_1>d_2,$ then $C$ is parametrized by 
$$ t\rightarrow \left(F_1^{\circ (d_1-d_2)}(t),\, (\alpha\circ\mu')(t)\right).$$ 
Since Corollary \ref{svin2} implies that $$F_1^{\circ (d_1-d_2)}\circ \mu'^{-1}\circ \alpha^{-1}=\mu''\circ F_1^{\circ (d_1-d_2)}\circ \alpha^{-1}$$ for some $\mu''\in \Aut(F_1^{\circ d})$, we see that   in this case $C$ is also  parametrized by 
 
\be \l{sl2}  t\rightarrow \left((\mu''\circ F_1^{\circ (d_1-d_2)}\circ\alpha^{-1})(t),\,t\right)\ee  for some $\mu''\in \Aut(F_1^{\circ d}).$  
This proves the ``only if'' part of the theorem.

In the other direction, let us assume  that \eqref{kosh1} holds and $C$ is a curve parametrized by 
$$t\rightarrow \left(t,\,(\alpha\circ \mu\circ F_1^{\circ s})(t)\right)$$ 
for some  $\mu\in \Aut(F_1^{\circ d}),$ M\"obius transformation $\alpha$,  and  $s\geq 0$. 
Since 
$$F_2^{\circ d}\circ (\alpha\circ \mu\circ F_1^{\circ s})=
\alpha\circ F_1^{\circ d}\circ \mu\circ F_1^{\circ s}=\alpha\circ \mu\circ F_1^{\circ d}\circ F_1^{\circ s}=(\alpha\circ \mu\circ F_1^{\circ s}) \circ  F_1^{\circ d},$$ in this case the 
 diagram 
\be  \l{kaba} 
\begin{CD} 
(\C\P^1)^2 @>(B,B)>>(\C\P^1)^2 \\ 
@V (X_1,X_2)  VV @VV  (X_1,X_2) V\\ 
 (\C\P^1)^2 @>(F_1^{\circ d},F_2^{\circ d})>> (\C\P^1)^2
\end{CD}
\ee
commutes for $$B=F_1^{\circ d},\ \ \ X_1=z,\ \ \ X_2=\alpha\circ \mu\circ F_1^{\circ s},$$ implying that \eqref{lass} holds. Similarly, if $C$ is parametrized by 
$$t\rightarrow \left((\mu\circ F_1^{\circ s}\circ\alpha^{-1})(t),\,t\right),$$
then it follows from 
$$F_1^{\circ d}\circ(\mu\circ F_1^{\circ s}\circ\alpha^{-1}) =\mu
\circ F_1^{\circ d}\circ F_1^{\circ s}\circ \alpha^{-1}=
\mu\circ F_1^{\circ s}\circ \alpha^{-1}\circ \alpha\circ  F_1^{\circ d}\circ \alpha^{-1}=$$ $$=(\mu\circ F_1^{\circ s}\circ\alpha^{-1})\circ F_2^{\circ d}$$ that diagram 
\eqref{kaba} commutes for 
 $$B=F_2^{\circ d},\ \ \ X_1=\mu\circ F_1^{\circ s}\circ\alpha^{-1},\ \ \ X_2=z. \eqno{\Box}$$

\noindent{\it Proof of Theorem \ref{t6}.} By Lemma \ref{lom1} and Lemma \ref{lom2},  there exists a Zariski open set $U$ in  ${\rm Rat}_m$ such that every $F\in U$  is simple and  the   group $G(F)$ is trivial. By Theorem \ref{t2}, this implies that for every $F\in U$  the   group $\Aut_{\infty}(F)$ is also trivial. Since for simple $F_1,$ $F_2$ equality \eqref{kosh1} yields equality \eqref{reby}, it follows now from Theorem \ref{t5} that if  $F_1,F_2\in U$, then $(F_1,F_2)$-periodic curves exist  if and  
only if  
\be \l{holds} F_2=\alpha \circ F_1\circ\alpha^{-1}\ee for some M\"obius transformation $\alpha$. 
Furthermore, these periodic curves have the form 
 $$y=(\alpha\circ F_1^{\circ s})(x), \ \ \ \ \ x=( F_1^{\circ s}\circ\alpha^{-1})(y),$$ and it is  easy to check arguing as above that  these  curves  are 
$(F_1,F_2)$-invariant. 
 \qed 
\section{Indecomposable functions with non-trivial decompositions of iterates}
Let $F$ be an indecomposable rational function. In  this section, we give a number of conditions implying that some iterate $F^{\circ k}$, $k>1$, of  $F$ has a decomposition into a composition of indecomposable rational functions that is not equivalent to $F^{\circ k}$. For brevity, in this case we will say that $F^{\circ k}$ has a {\it non-trivial decomposition}.  We also give explicit examples of simple rational functions of degrees 2 and 3 for which  Theorems \ref{t1} - \ref{t2} and Theorems \ref{t4} - \ref{t5} do not hold.  

Let us recall that for a rational function  $F$ the group 
$ G(F)$ is defined as the group  of all M\"obius transformations $\sigma$ such that $$ F\circ \sigma=\nu\circ F$$ for some  M\"obius transformations $\nu$. 
 Along with the group $G(F)$ we will consider its subgroup $\Sigma(A)$ consisting of all M\"obius transformations $\sigma$ such that $ A\circ \sigma= A$. 
We recall that for  a finite subgroup $G$ of $\Aut(\C\P^1)$ an {\it invariant function}  for $G$ is a rational function $\theta_G$ such that the equality $$\theta_G(x)=\theta_G(y), \ \ \ \ \ x,y\in \C\P^1,$$  holds if and only if there exists $\sigma\in G$ such that $\sigma(x)=y.$ 
Such a function always exists and is defined in a unique way up to the transformation $\theta\rightarrow \mu \circ \theta,$ where $\mu$ is a M\"obius transformation. 
Obviously,  the degree of $\theta_G$  is equal to the order of $G$, and  
it follows easily from  
 the L\"uroth theorem that a rational function   $F$ is a rational function in $\theta_G$
if and only if $G\subseteq \Sigma(F)$.

\bt \l{tt1} 
Let $F$ be an indecomposable rational function such that the group $G(F^{\circ k_0})$ for some $k_0>1$ 
contains an element that does not  belong to the group $G(F)$. Then the iterate 
$F^{\circ k_0}$ has a non-trivial decomposition.
\et 
\pr Let $\alpha$ be a M\"obius transformation such that  $\alpha\in G(F^{\circ k_0})$ but $\alpha\not \in G(F)$. Then 
$$F^{\circ k_0}\circ \alpha=\nu \circ F^{\circ k_0}$$ for some  M\"obius transformation $\nu$, implying that 
$$F^{\circ k_0}=(\nu^{-1}\circ F)\circ F \circ \dots \circ F \circ (F\circ \alpha).$$ Moreover, the decomposition on the right side of this equality   
is non-trivial since otherwise $F=\mu\circ (F\circ \alpha)$ for some   M\"obius transformation  $\mu$,  in contradiction with the assumption $\alpha \not \in G(F)$. \qed

\bt\l{tt2}
Let $F$ be an indecomposable rational function of degree $n\geq 2$ that is  not conjugate to $z^{\pm n}$ such that the group $\Sigma(F^{\circ k_0})$ for some $k_0>1$ contains an element that does not  belong to the groups $\Sigma(F^{\circ k})$, $1\leq k< k_0$.  Then the iterate 
$F^{\circ k_0}$ has a non-trivial decomposition.
\et   
\pr Let $\sigma$  be such an element. 
Since $\langle\sigma\rangle$ is a subgroup of $\Sigma(F^{\circ k_0})$, there exists a rational function $R$ such that 
$F^{\circ k_0}=R\circ \theta_{\langle\sigma\rangle}$. Assume for a moment that   $$\vert \langle\sigma\rangle \vert <(\deg F)^{k_0}.$$  
Since $\deg \theta_{\langle\sigma\rangle}=\vert \langle\sigma\rangle \vert,$ in this case $
\deg R\geq 2$. Let now $$R=G_1\circ G_2\circ \dots \circ G_l \ \ \ {\rm and} \ \ \ \theta_{\langle\sigma\rangle}=H_1\circ H_2\circ \dots \circ H_t$$ 
 be any decompositions 
into compositions of indecomposable rational functions. Concatenating them  we 
obtain a decomposition 
\be \l{rs} F^{\circ k_0}=G_1\circ G_2\circ \dots \circ G_l\circ H_1\circ H_2\circ \dots \circ H_t\ee  of $ F^{\circ k_0}.$ 
If this decomposition is equivalent to $F^{\circ k_0}$, then 
$$F^{\circ k}=\mu \circ H_1\circ H_2\circ \dots \circ H_t=\mu\circ \theta_{\langle\sigma\rangle},$$ for some M\"obius transformation $\m$ and $k\geq 1$. Moreover, $k<k_0$ since $\deg R\geq 2.$ 
 Thus,  $\sigma\in \Sigma(F^{\circ k})$, where $k<k_0$, in contradiction with the condition. Hence, decomposition  \eqref{rs} is non-trivial. 

To finish the proof, we only must show that the equality  
\be \l{ifh} \vert \langle\sigma\rangle \vert =(\deg F)^{k_0}\ee for $\sigma\in \Aut(\C\P^1)$ with $\vert  \langle\sigma\rangle \vert<\infty$ 
 implies that $F$ is conjugate to $z^{\pm n}$.   
Since 
$$\theta_{\langle\sigma\rangle}=\alpha\circ z^{\vert \langle\sigma\rangle \vert }\circ \beta$$ 
 for some M\"obius transformations $\alpha$ and $\beta$, this reduces to showing that if 
\be \l{ies} F^{\circ k_0}= z^{n^{k_0}}\circ \beta\ee for some $k_0>1$ and M\"obius transformations $\beta$, then $F$ is conjugate to $z^{\pm n}$.

Clearly, if \eqref{ies} holds, then the preimage 
$\left(F^{\circ k_0}\right)^{-1}\{0,\infty\}$ contains only two points, implying that all the preimages $\left(F^{\circ k}\right)^{-1}\{0,\infty\}$, $1\leq k <k_0$,   also 
 contain only two points. Set $\{a,b\}=F^{-1}\{0,\infty\}$. Then it follows from \eqref{ies} that  
\be \l{th} \left|\left(F^{\circ (k_0-1)}\right)^{-1}\{0,\infty\}\right|=\left|\left(F^{\circ (k_0-1)}\right)^{-1}\{a,b\}\right|=2,\ee implying that  
\be \l{iim} \left| \left(F^{\circ (k_0-1)}\right)^{-1}\{0,\infty,a,b\}\right| \leq 4.\ee 

Let us recall now that the Riemann-Hurwitz formula implies that the preimage of a  finite set $S$ under a rational function $H$ of degree $d\geq 2$ contains at least $d(\vert S\vert -2)+2$ points,  and the equality is attained if and only if the set of critical values $c(H)$ of $H$ belongs to $S$. In particular, 
if $S$ contains at least three points, then $$\vert H^{-1}\{S\}\vert \geq d+2.$$ 
Thus,  if $n> 2$ or $k_0>2$, then 
\eqref{iim} is possible only if $\{0,\infty\}=\{a,b\}$, implying that $F$ is conjugate to $z^{\pm n}.$ On the other hand, if $n=k_0=2$, then it follows from \eqref{th} that $c(F)=\{0,\infty\}$ and $c(F)=\{a,b\}$, 
 implying again that $\{0,\infty\}=\{a,b\}$. \qed

As an example illustrating the above theorems, let us consider the function 
\be \l{inf} F=\frac{z^2-1}{z^2+1}.\ee Since $F=\frac{z-1}{z+1}\circ z^2$, the group $G(F)$ 
 consists of the transformations $c z^{\pm 1}$, $c\in \C\setminus\{0\}$ (see e. g. Section 2 of \cite{sym}). On the other hand, for 
$$F^{\circ 2}=-\,{\frac {2{z}^{2}}{{z}^{4}+1}},$$ the corresponding group 
$G(F^{\circ 2})$ contains the transformation $\mu=\frac{z+i}{z-i}$ satisfying $$F^{\circ 2}\circ \mu=\nu\circ F^{\circ 2}$$ for $\nu={\frac {-z+1}{-3\,z-1}}$. Hence, 
\be \l{d1} \frac{z^2-1}{z^2+1}\circ \frac{z^2-1}{z^2+1}=\left(\nu^{-1}\circ\frac{z^2-1}{z^2+1}\right)\circ \left(\frac{z^2-1}{z^2+1}\circ \mu\right)=-{\frac {{z}^{2}}{{z}^{2}-2}}\circ {\frac {2\,iz}{{z}^{2}-1}},\ee where the decomposition on the right side  
is non-trivial by Theorem \ref{tt1}. 

Furthermore, the transformation $\delta: z\rightarrow \frac{1}{z}$ obviously satisfies 
$$F\circ \delta =- F, \ \ \ F^{\circ 2}\circ \delta = F.$$ Thus, $\delta\in \Sigma(F^{\circ 2})$ but $\delta\not\in \Sigma(F)$,   and therefore  
$F^{\circ 2}$ is a rational function in $$\theta_{\langle \delta\rangle}=z+\frac{1}{z}.$$ The corresponding non-trivial decomposition provided by Theorem \ref{tt2} is given by the formula 
\be \l{d2} \frac{z^2-1}{z^2+1}\circ \frac{z^2-1}{z^2+1}=-\frac{2}{{z}^{2}-2} \circ \left(z+\frac{1}{z}\right).\ee  Notice that although $\delta$ does not belong to $\Sigma(F)$, it belongs to $G(F)$. Thus, 
we cannot apply  to $\delta$ previous Theorem \ref{tt1}. 

Let us remark that decompositions of $F^{\circ 2}$ on the right sides of \eqref{d1} and \eqref{d2} are not equivalent to each other. Indeed, for equivalent decompositions \eqref{lll}  the sets of critical points of $F_1$ and $G_1$ are equal. On the other hand, the sets of critical points of 
the  functions $\frac {2\,iz}{{z}^{2}-1}$ and $z+\frac{1}{z}$  are $\{-i,i\}$ and $\{-1,1\}$. Notice that since the set of critical points of $F$ is $\{0,\infty\}$, this also gives another proof of the fact that decompositions \eqref{d1} and \eqref{d2} are non-trivial. 
Finally, let us mention that $F^{\circ 2}$ has no other non-trivial decompositions. To see this, let us observe that  the function $F^{\circ 2}$ is an invariant function for  
the Klein four group $V\cong \Z/2\Z\times \Z/2\Z$ generated by the transformations 
$z\rightarrow \frac{1}{z}$ and $z\rightarrow -z.$
Since  $\Mon(\theta_G)\cong G$ for any finite subgroup $G$  of $\Aut(\C\P^1)$, this implies  that $\Mon(F^{\circ 2})\cong V$. Therefore, as the group $V$ has three proper imprimitivity systems,  
 $F^{\circ 2}$ has three non-equivalent decompositions (one of which corresponds to the decomposition 
$F^{\circ 2}$ itself). 

In relation with Theorem \ref{tt1} and Theorem \ref{tt2}, let us mention that according to the recent results of \cite{sym} for any rational function $A$ of degree $n\geq 2$ that is not conjugate to $z^n,$ the orders of the groups 
$G(A^{\circ k})$, $k\geq 2,$  are finite and uniformly bounded
in terms of $n$ only. 

To construct further examples of rational functions whose iterates admit non-trivial decompositions, we will use Latt\`es maps. More precisely, we will use those Latt\`es maps that can be obtained as projections on the $x$-coordinate of self-isogenies of elliptic curves. 
We emphasize that this class of Latt\`es maps does not exhaust the entire class Latt\`es maps as it was defined in Section \ref{44}. However, for brevity, in the rest of the article we will call by  Latt\`es maps only these particular  Latt\`es maps, which are defined in detail below. 
Notice that this time we use another, more common definition of Latt\`es maps. 

Let $\f C$ and $\t{\f C}$ be elliptic curves over $\C$ defined in the short Weierstrass form. We recall that 
an {\it isogeny} between $\f C$ and $\t{\f C}$ is a morphism $\psi:\f C\rightarrow \t{\f C}$ that sends 
the identity element of the  group  $\f C$ to the  identity element of the  group  $\t{\f C}$. 
 Such a morphism is necessarily  a homomorphism of groups (see e. g.  \cite{sil}). 
Thus, $\psi(-x)=-\psi(x)$,  implying that 
there exists a rational function $F$ such that 
the diagram
\be \l{dig1} 
\begin{CD}
\f C @>\psi>> \t{\f C} \\
@VV x  V @VV x V\\ 
\C\P^1 @>F >> \ \ \C \P^1\, 
\end{CD}
\ee
commutes. 
If $\f C=\t{\f C}$ and $\psi$ is an endomorphism, we will call the corresponding rational function $F$ a Latt\`es map. 
For example, for any elliptic curve $\f C$, the multiplication by $n$ on $\C$ induces 
an endomorphism $[n]:\f C\rightarrow \f C$ of degree $n^2$. We will denote by $F_{\f C,n}$
the corresponding Latt\`es map of degree $n^2$,  
which makes the diagram 
\be \l{dig2} 
\begin{CD}
\f C @>[n]>> \f C \\
@VV x  V @VV x V\\ 
\C\P^1 @>F_{\f C,n} >> \ \ \C \P^1\, 
\end{CD}
\ee
commutative.

Below, we will use the following results about isogenies (see e.g. \cite{sil}, Ch. III).
Let $\psi: \f C\rightarrow \t{\f C}$ be a non-zero isogeny of degree $n$. Then its kernel  $\Gamma$   is a subgroup of order $n$ in  
$\f C$. Moreover,  for any  subgroup of  $\f C$  there exists an isogeny  $\psi:\f C\rightarrow \tt {\f C}$ such that $\ker \psi=\Gamma$, and this isogeny is defined in a unique way up to an isomorphism of $\tt {\f C}$.   Finally, for any isogeny 
$\psi:\f C\rightarrow \t{\f C}$ of degree $n$ there exists 
 a unique  {\it dual} isogeny  $\widehat\psi:\, \t{\f C}\rightarrow \f C$ such that $\widehat \psi \circ \psi =[n]$  on $\f C$ and $ \psi \circ \widehat \psi =[n]$  on $\tt{ \f C}$.

\bt\l{tla}
Let $p$ be a prime number, and $\f C$ an elliptic curve such that the multiplication by $i\sqrt{p}$ on $\C$ induces  an endomorphism of $\f C.$ 
Then the corresponding  Latt\`es map $F$  is indecomposable, and the iterate 
$F^{\circ 2}$ has a non-trivial decomposition.
\et   
\pr Clearly,  $F^{\circ 2}$ 
 makes the diagram 
\be  
\begin{CD}
\f C @>-[p]>> \f C \\
@VV x  V @VV x V\\ 
\C\P^1 @>F^{\circ 2} >> \ \ \C \P^1\, 
\end{CD}
\ee
commutative.  Since the change of the sign of $\psi$ in \eqref{dig1} obviously does not affect the corresponding Latt\`es map, this implies that $F^{\circ 2}=F_{\f C,p}$. Therefore, as any rational function of prime degree is clearly indecomposable, to prove the theorem it is enough to show that the function $F_{\f C,p}$ has more than one equivalence class of decompositions into compositions of indecomposable rational functions. 

Let us show that starting from any finite subgroup $\Gamma$ of $\f C$ of order $p$  one can construct a decomposition of $F_{\f C,p}$ into a composition of rational functions of degree $p$. Let $\psi:\f C\rightarrow \f C_{\Gamma}$ be an isogeny such that $\ker \psi=\Gamma$. Then 
there exists a rational function $V_{\Gamma}$ such that 
the diagram
$$
\begin{CD}
\f C @>\psi_{\Gamma}>> {\f C_{\Gamma}} \\
@VV x V @VV x V\\ 
\C\P^1 @>V_{\Gamma} >> \ \ \C \P^1\, 
\end{CD}
$$
commutes. Similarly, for the dual isogeny $\widehat\psi_{\Gamma}:\, {\f C_{\Gamma}}\rightarrow \f C$  there exists   
a rational function $U_{\Gamma}$ such that 
 the diagram 
$$
\begin{CD}
\f C_{\Gamma} @>\widehat \psi_{\Gamma}>> \f C \\
@VV x V @VV x  V\\ 
\C\P^1 @>U_{\Gamma} >> \ \ \C \P^1\, 
\end{CD}
$$ commutes. Gluing  these diagrams we obtain a decomposition 
\be \l{gam} F_{\f C,p}=U_{\Gamma}\circ V_{\Gamma}.\ee  

Let us prove now that if $\Gamma_1\neq \Gamma_2$, then the decompositions 
$U_{\Gamma_1}\circ V_{\Gamma_1}$ and $ U_{\Gamma_2}\circ V_{\Gamma_2}$ are not equivalent. 
Let us consider the maps  $V_{\Gamma_1}\circ U_{\Gamma_1}$ and $ V_{\Gamma_2}\circ U_{\Gamma_2}$, which make the diagrams 
$$
\begin{CD}
\f C_{\Gamma_1} @>\psi_{\Gamma_1}\circ \widehat \psi_{\Gamma_1}>> \f C_{\Gamma_1} \\
@VV x V @VV x  V\\ 
\C\P^1 @>V_{\Gamma_1}\circ U_{\Gamma_1} >> \ \ \C \P^1\,, 
\end{CD} 
\ \ \ \ \ \ \ \ \ \ \ 
\begin{CD}
\f C_{\Gamma2} @> \psi_{\Gamma_2}\circ \widehat \psi_{\Gamma_2}>> \f C_{\Gamma_2} \\
@VV x V @VV x  V\\ 
\C\P^1 @>V_{\Gamma_2}\circ U_{\Gamma_2} >> \ \ \C \P^1\, 
\end{CD}
$$
commutative. Clearly,  $V_{\Gamma_1}\circ U_{\Gamma_1}$ and $ V_{\Gamma_2}\circ U_{\Gamma_2}$ are Latt\`es maps. Moreover, if the decompositions $U_{\Gamma_1}\circ V_{\Gamma_1}$ and $ U_{\Gamma_2}\circ V_{\Gamma_2}$ are equivalent, then these Latt\`es maps are conjugate. 
However, for conjugate Latt\`es maps the corresponding elliptic curves are isomorphic 
(see e.g. \cite{sildyn}, Theorem 6.46). Thus, if the above decompositions are equivalent, then  $\f C_{\Gamma_1}\cong \f C_{\Gamma_2},$  implying that $\Gamma_1= \Gamma_2$. 

To finish the proof, we only must show that there exist at least two different subgroups of order $p$ in $\f C$. In fact, it is easy to see that there exist exactly $p+1$ such subgroups. Indeed, any  subgroup of order $p$ is cyclic and is contained in the subgroup of points of $\f C$ whose order divides $p$, that is, in the kernel of $[p]$.  Thus, the number of subgroups of order $p$  is equal to the number of cyclic subgroups of order $p$ in  $\Z/p\Z\times \Z/p\Z$. In turn, this number is equal to the number  of elements of order $p$ in $\Z/p\Z\times \Z/p\Z$, which is equal to $p^2-1$, divided by the number of elements generating the same subgroup, which is equal to $p-1$. \qed

Notice that Latt\`es maps satisfying conditions of Theorem \ref{tla} exist for every prime $p$. To see this, it is enough to observe that $i\sqrt{p}$ is an endomorphism of an elliptic curve corresponding to the lattice generated by $1$ and $i\sqrt{p}$.   

 Examples illustrating Theorem \ref{tla}  for $p=2$ and $p=3$ are given by the functions 
\be \l{inl} L={\frac {2\,\sqrt {2}x-{x}^{2}-1}{2x}} \ \ \ \ \ {\rm and} \ \ \ \ \ P={\frac {6x}{{x}^{3}-2}}\ee 
(see \cite{mil2}). 
Using  V\'elu's formulas for isogenies (\cite{velu})  
 or a brute force calculation, one can find the following non-trivial decompositions of 
their second iterates:
\be \l{deca1} L^{\circ 2}=U\circ V,\ee 
where 
$$U=\frac{{x}^{2}}{4( x- \frac{1}{\sqrt {2}-1})}, \ \ \  V=\frac{x^{2}-1}{x-\sqrt {2}+1},$$ 
and \be \l{deca2} P^{\circ 2}=Q\circ R,\ee where 
$$Q=-{\frac {23328x}{{x}^{3}+216\,\sqrt [3]{2}{x}^{2}+3888\,{2}^{2/3}x-
93312}}, \ \ \ R={\frac {36x \left( {2}^{2/3}{x}^{2}-4\,x+2\,\sqrt [3]{2} \right) }{{
2}^{2/3}{x}^{2}+2\,x+2\,\sqrt [3]{2}}}.$$
As above, to see directly that these decompositions are indeed non-trivial it is enough to compare 
the sets of critical points of $L$ and  $V$ for \eqref{deca1} and the sets of critical points of $P$ and $R$ for \eqref{deca2}.  In the first case, one can check that 
 the corresponding sets are 
$$\{-1,1\} \ \ \ \ {\rm and} \ \ \ \   \{\sqrt {2}-1+i\sqrt {-2+2\,\sqrt {2}},\sqrt {2}-1-i\sqrt {-2+2\,\sqrt {2}}\}.$$ In the second case, it is enough to observe that $\infty$ is a critical point of $P$, but is not a critical point of $R$. 

It is clear that any rational function of degree 2 is simple. Moreover, one can check that the function  $P$ 
in \eqref{inl} has four critical values (in fact it is true for any Latt\`es map) and hence  is  simple by Lemma \ref{bl}. Thus, the functions given by \eqref{inf} and \eqref{inl} provide us with 
counterexamples to Theorem \ref{t1} for $m$ equal 2 and 3. Moreover,  for these values of $m$, functions given by  \eqref{inl} give counterexamples also  to Theorem \ref{t2}. Indeed, for a Latt\`es map $F$ the semigroups $E(F)$ and $C_{\infty}(F)$ are not finitely generated (see \cite{mil2}). On the other hand, the  group 
$\Aut_{\infty}(F)$ is finite for any rational function $F$ (see \cite{sym}). 

A simple counterexample to Theorem \ref{t4} 
 for $m=2$ is obtained from 
the semiconjugacy 
$$
\begin{CD}
\C\P^1 @>z^m>> \C\P^1 \\
@V \frac{1}{2}\left(z+\frac{1}{z}\right) VV @VV \frac{1}{2}\left(z+\frac{1}{z}\right)  V\\ 
\f \C\P^1 @>T_m>> \f\C\P^1\, 
\end{CD}
$$
 for $m=2$. Indeed, it is clear that $z^m$ and $T_m$ are not conjugate. On the other hand, $T_2$ is simple. Furthermore, the commutative diagram  
$$
\begin{CD} 
(\C\P^1)^2 @>(z^2,z^2)>>(\C\P^1)^2 \\ 
@V \left(z,\frac{1}{2}\left(z+\frac{1}{z}\right)\right)  VV @VV \left(z,\frac{1}{2}\left(z+\frac{1}{z}\right)\right)  V\\ 
 (\C\P^1)^2 @>(z^2,T_2)>> (\C\P^1)^2\,
\end{CD}
$$
 defines for non-conjugate simple rational functions $z^2$ and $T_2$ an irreducible $(z^2,T_2)$-invariant curve parametrized by $z\rightarrow \left(z,\frac{1}{2}\left(z+\frac{1}{z}\right)\right)$, providing a counterexample to  
 Theorem \ref{t5}. 

Counterexamples to Theorem \ref{t4} and Theorem \ref{t5} 
 for $m=3$ can be obtained from 
the semiconjugacy  
$$  
\begin{CD}
\C\P^1 @>z\left(\frac{z^2-a}{z^2-b}\right) >> \C\P^1\\
@VV {z^2} V @VV z^2 V\\ 
\C\P^1 @>z\left(\frac{z-a}{z-b}\right)^2 >> \C\P^1,
\end{CD}
$$
where $a,b\in \C$, which is a particular case of the semiconjugacy 
\be 
\begin{CD}
\C\P^1 @>z^rR(z^n) >> \C\P^1\\
@VV {z^n} V @VV z^n V\\ 
\C\P^1 @>z^rR^n(z) >> \C\P^1,
\end{CD}
\ee
where $R$ is an arbitrary rational function and $r,n$ are 
integer positive numbers. 

Setting, for example, $a=2,b=3$ and observing that $$A=z\left(\frac{z-2}{z-3}\right)^2$$ has three fixed points 
$0,\infty, 5/2$, while $$B=z\left(\frac{z^2-2}{z^2-3}\right)$$ has only two fixed points $0$ and $\infty$,   we see that $A$ and $B$ are not conjugate. On the other hand, since $A$ has four critical values $0,\infty, 32/3,1/4$, it is a simple rational function. Thus, we obtain 
 a counterexample to  Theorem \ref{t4}. A counterexample to  Theorem \ref{t5} is obtained from the diagram 
$$
\begin{CD} 
(\C\P^1)^2 @>(B,B)>>(\C\P^1)^2 \\ 
@V (z,z^2)  VV @VV (z,z^2)  V\\ 
 (\C\P^1)^2 @>(B,A)>> (\C\P^1)^2\,
\end{CD}
$$
in the same way as above, taking into account that $B$ is also simple since it has four critical values $\pm \frac{1}{2},$ $\pm \frac{4\sqrt{6}}{3}.$  

In conclusion, we mention that if $F$ is {\it decomposable}, then the problem of describing the whole totality of its iterates seems to be even more complicated than in the indecomposable case, and a ``qualitative'' description comes to the fore. An example of such a description is the result of   (\cite{mz}), which states that if $F$ is a {\it polynomial}  of degree $n\geq 2$ not conjugate to $z^n$ or to $\pm T_n$, then decompositions of its iterates can be obtained from decompositions of a single iterate 
 $F^{\circ N}$ for $N$ big enough in the following sense: 
 there exists an integer $N\geq 1 $ such that any decomposition of  $F^{\circ d}$ with $d\geq N$
has the form 
$$X=F^{\circ k_1}\circ X', \ \ \ \ \  Y=Y'\circ F^{\circ k_2},$$ 
where $F^{\circ N}= X'\circ Y'$ and $k_1,k_2\geq 0$ (see \cite{mz} and also \cite{pj}, \cite{tame} for different proofs of  this fact).  

Conjecturally, the result of \cite{mz} remains true for any  non-special rational function $F$. To the date, this conjecture  is proved for ``tame'' rational functions, that is, for  
the functions $F$ satisfying the following condition: the algebraic curve $$F(x)-F(y)=0$$ has no factors of genus zero or one distinct from the diagonal (\cite{tame}).  
Since any simple rational function $F$ of degree $m\geq 4$ is tame by Theorem \ref{goo1}  and formula \eqref{burry},  
this also shows  that the conjecture is true for simple and general rational functions of degree at least four. On the other hand, Theorem \ref{t1} shows that for such  functions  it is true already for $N=1.$

\end{document}